\documentclass[10pt]{article}

\usepackage[dvips]{epsfig} 
\usepackage{amssymb,amsmath,amsthm,mathrsfs} 
\usepackage{graphicx}

\usepackage[authoryear,sort&compress]{natbib}
\usepackage{url}
\usepackage{enumerate}
\usepackage{colordvi,color,pspicture}
\usepackage{psfrag}

\usepackage{setspace}

\graphicspath{{Figures/}}

\newcommand{\Y}{\mathcal{Y}}
\newcommand{\X}{\mathcal{X}}
\newcommand{\U}{\mathcal{U}}
\newcommand{\I}{\mathbf{I}}
\newcommand{\C}{\mathcal{C}}
\newcommand{\mI}{\mathcal{I}}
\newcommand{\y}{\mathsf{y}}

\newcommand{\F}{\mathcal F}
\newcommand{\g}{\gamma}
\newcommand{\G}{\Gamma}
\newcommand{\N}{\mathcal{N}}
\newcommand{\V}{\mathcal{V}}
\newcommand{\A}{\mathcal{A}}
\newcommand{\R}{\mathbb{R}}
\newcommand{\mS}{\mathcal{S}}
\newcommand{\E}{\mathbf{E}}

\newcommand{\D}{\mathscr D}
\newcommand{\ve}{\varepsilon}

\newcommand{\al}{\alpha}
\newcommand{\be}{\beta}

\newcommand{\TY}{T(\Y_3)}

\DeclareMathOperator{\argmin}{argmin}

\DeclareMathOperator{\BER}{BER}
\DeclareMathOperator{\BIN}{BIN}

\theoremstyle{plain}
\newtheorem{theorem}{Theorem}[section]
\newtheorem{lemma}[theorem]{Lemma}
\newtheorem{proposition}[theorem]{Proposition}
\newtheorem{corollary}[theorem]{Corollary}

\theoremstyle{definition}

\newtheorem{example}[theorem]{Example}

\theoremstyle{remark}

\newtheorem{remark}[theorem]{Remark}

\begin{document}

\title{
Technical Report \# KU-EC-10-2\\
The Distribution of the Domination Number of a Family of Random Interval Catch Digraphs}
\author{
Elvan Ceyhan\thanks{Address:
Department of Mathematics, Ko\c{c} University, 34450 Sar{\i}yer, Istanbul, Turkey.
e-mail: elceyhan@ku.edu.tr, tel:+90 (212) 338-1845, fax: +90 (212) 338-1559.
}
}

\date{\today}

\maketitle

\begin{abstract}
\noindent
We study a new kind of proximity graphs called
proportional-edge proximity catch digraphs (PCDs)
in a randomized setting.
PCDs are a special kind of random catch digraphs that have been developed recently
and have applications in statistical pattern classification and spatial point pattern analysis.
PCDs are also a special type of intersection digraphs;
and for one-dimensional data,
the proportional-edge PCD family is also a family of random interval catch digraphs.
We present the exact (and asymptotic) distribution of the domination number of
this PCD family for uniform (and non-uniform) data in one dimension.
We also provide several extensions of this random catch digraph by relaxing
the expansion and centrality parameters, thereby determine the parameters
for which the asymptotic distribution is non-degenerate.
We observe sudden jumps
(from degeneracy to non-degeneracy or
from a non-degenerate distribution to another)
in the asymptotic distribution of the domination number at certain parameter values.
\end{abstract}

\vspace{0.1 in}



\newpage
\section{Introduction}
\label{sec:intro}
The proximity catch digraphs (PCDs) were motivated by their applications
in pattern classification and spatial pattern analysis,
hence they have become focus of considerable attention recently.
The PCDs are a special type of proximity graphs which were introduced by \cite{toussaint:1980}.
A \emph{digraph} is a directed graph with
vertex set $V$ and arcs (directed edges) each of which is from one
vertex to another based on a binary relation.
Then the pair $(p,q) \in V \times V$ is an ordered pair
which stands for an \emph{arc} from vertex $p$ to vertex $q$ in $V$.
For example, \emph{nearest neighbor (di)graph} which is defined by placing
an arc between each vertex and its nearest neighbor is a proximity digraph
where vertices represent points in some metric space (\cite{paterson:1992}).
PCDs are \emph{data-random digraphs} in which
each vertex corresponds to a data point
and arcs are defined in terms of some bivariate relation on the data.

The PCDs are closely related to the class cover problem of \cite{cannon:2000}.
Let $(\Omega,\mathcal M)$ be a measurable space and $\X_n=\{X_1,X_2,\ldots,X_n\}$ and
$\Y_m=\{Y_1,Y_2,\ldots,Y_m\}$ be two sets of $\Omega$-valued random variables
from classes $\X$ and $\Y$, respectively, with joint probability distribution $F_{X,Y}$.
Let $d(\cdot,\cdot):\Omega\times \Omega \rightarrow [0,\infty)$ be any distance function.
The \emph{class cover problem} for a target class, say $\X$, refers to finding a collection of neighborhoods,
$N_i$ around $X_i$ such that
(i) $\X_n \subseteq  \bigl(\cup_i N_i \bigr)$ and (ii) $\Y_m \cap \bigl(\cup_i N_i \bigr)=\emptyset$.
A collection of neighborhoods satisfying both conditions is called a {\em class cover}.
A cover satisfying (i) is a {\em proper cover} of $\X_n$
while a cover satisfying (ii) is a {\em pure cover} relative to $\Y_m$.
This article is on the {\em cardinality of smallest class covers}; that is,
class covers satisfying both (i) and (ii) with the smallest number of neighborhoods.
See \cite{cannon:2000} and \cite{priebe:2001} for more on the class cover problem.

The first type of PCD was class cover catch digraph (CCCD) introduced by \cite{priebe:2001}
who gave the exact distribution of its domination number for
uniform data from two classes in $\R$.
\cite{devinney:2002a}, \cite{marchette:2003}, \cite{priebe:2003b}, \cite{priebe:2003a},
\cite{devinney:2006} extended the CCCDs to higher dimensions and demonstrated that CCCDs
are a competitive alternative to the existing methods in classification.
Furthermore,
\cite{devinney:2002b} proved a SLLN result for the one-dimensional class cover problem;
\cite{wiermanSLLN:2008} provided a generalized SLLN result
and \cite{xiangCLT:2009} provided a CLT result for
CCCD based on one-dimensional data.
However, CCCDs have some disadvantages in higher dimensions;
namely, finding the minimum dominating set for CCCDs is an NP-hard problem in general,
although a simple linear time algorithm is available
for one dimensional data (\cite{priebe:2001});
and the exact and the asymptotic distributions of the domination number
of the CCCDs are not analytically tractable in multiple dimensions.
\cite{ceyhan:CS-JSM-2003,ceyhan:dom-num-NPE-SPL} introduced the central similarity proximity maps
and proportional-edge proximity maps for data in $\R^d$ with $d>1$
and the associated random PCDs
with the purpose of avoiding the above mentioned problems.
The asymptotic distribution of the domination number of the proportional-edge PCD
is calculated for data in $\R^2$ and
then the domination number is used as a statistic
for testing bivariate spatial patterns
(\cite{ceyhan:dom-num-NPE-SPL}, \cite{ceyhan:dom-num-NPE-Spat2010}).
The relative density of these two PCD families is also calculated
and used for the same purpose
(\cite{ceyhan:arc-density-PE} and \cite{ceyhan:arc-density-CS}).
Moreover, the distribution of the domination number
of CCCDs is derived for non-uniform data (\cite{ceyhan:dom-num-CCCD-NonUnif}).

In this article,
we provide the exact (and asymptotic) distribution of the domination number of
proportional-edge PCDs for uniform (and non-uniform) one-dimensional data.
First, some special cases and bounds for the domination number of proportional-edge PCDs
is presented, then the domination number is investigated
for uniform data in one interval (in $\R$) and the analysis is generalized to
uniform data in multiple intervals and to non-uniform data in one and multiple intervals.
These results can be seen as generalizations of the results of \cite{ceyhan:dom-num-CCCD-NonUnif}.
Some trivial proofs are omitted,
shorter proofs are given in the main body of the article;
while longer proofs are deferred to the Appendix.

We define the proportional-edge PCDs
and their domination number in Section \ref{sec:prop-edge-PCD},
provide the exact and asymptotic distributions of the domination number of proportional-edge PCDs
for uniform data in one interval in Section \ref{sec:gamma-dist-uniform},
discuss the distribution of the domination number for
data from a general distribution in Section \ref{sec:non-uniform}.
We extend these results to multiple intervals in Section \ref{sec:dist-multiple-intervals},
and provide discussion and conclusions in Section \ref{sec:disc-conclusions}.
For convenience in notation and presentation,
we resort to non-standard extended (perhaps abused) forms of Bernoulli and Binomial distributions,
denoted $\BER(p)$ and $\BIN(n,p)$, respectively,
where $p$ is the probability of success and $n$ is the number of trials.
Throughout the article,
we take $p \in [0,1]$ (unlike $p \in (0,1)$)
and if $X \sim \BER(p)$, then $P(X=1)=p$ and $P(X=0)=1-p$.
If $Y \sim \BIN(n,p)$,
then
$P(Y=k)=(n \atop k)p^k(1-p)^{n-k}$ for $p \in (0,1)$ and $k \in \{0,1,2,\ldots,n\}$
and
$P(Y=n)=1$ for $p=1$ and $P(Y=0)=1$ for $p=0$.

\section{Proportional-Edge Proximity Catch Digraphs}
\label{sec:prop-edge-PCD}
Consider the map $N:\Omega \rightarrow \wp(\Omega)$
where $\wp(\Omega)$ represents the power set of $\Omega$.
Then given $\Y_m \subseteq \Omega$,
the {\em proximity map}
$N(\cdot): \Omega \rightarrow \wp(\Omega)$
associates with each point $x \in \Omega$
a {\em proximity region} $N(x) \subseteq \Omega$.
For $B \subseteq \Omega$, the $\G_1$-region is the image of the map
$\G_1(\cdot,N):\wp(\Omega) \rightarrow \wp(\Omega)$
that associates the region $\G_1(B,N):=\{z \in \Omega: B \subseteq  N(z)\}$
with the set $B$.
For a point $x \in \Omega$, we denote $\G_1(\{x\},N)$ as $\G_1(x,N)$.
Notice that while the proximity region is defined for one point,
a $\G_1$-region is defined for a set of points.
The {\em data-random PCD} has the vertex set $\V=\X_n$
and arc set $\A$ defined by $(X_i,X_j) \in \A$ iff $X_j \in N(X_i)$.

Let $\Omega=\R$ and $Y_{(i)}$ be the $i^{th}$ order statistic
(i.e., $i^{th}$ smallest value) of $\Y_m$ for $i=1,2,\ldots,m$
with the additional notation for $i \in \{0,m+1\}$ as
$$-\infty =: Y_{(0)}<Y_{(1)}< \ldots <Y_{(m)}< Y_{(m+1)}:=\infty.$$
Then $Y_{(i)}$ partition $\R$ into $(m+1)$ intervals
which is called the \emph{intervalization} of $\R$ by $\Y_m$.
Let also that
$\mI_i:=\left( Y_{(i)},Y_{(i+1)} \right)$ for $i \in \{0,1,2,\ldots,m\}$
and $M_{c,i}:=Y_{(i)}+c(Y_{(i+1)}-Y_{(i)})$
(i.e., $M_{c,i} \in \mI_i$ such that $c \times 100$ \% of
length of $\mI_i$ is to the left of $M_{c,i}$).
We define the proportional-edge proximity region
with the expansion parameter $r \ge 1$ and
centrality parameter $c \in [0,1]$
for two one-dimensional data sets, $\X_n$ and $\Y_m$,
from classes $\X$ and $\Y$, respectively, as follows.
For $x \in \mI_i$ with $i \in \{1,2,\ldots,m-1\}$
\begin{equation}
\label{eqn:NPEr-general-defn1}
N(x,r,c)=
\begin{cases}
\left( Y_{(i)}, Y_{(i)}+r\,\left( x-Y_{(i)} \right) \right) \cap \mI_i      & \text{if $x \in (Y_{(i)},M_{c,i})$,}\\
\left( Y_{(i+1)}-r\left(Y_{(i+1)}-x\right), Y_{(i+1)} \right)  \cap \mI_i   & \text{if $x \in \left( M_{c,i},Y_{(i+1)} \right)$,}
\end{cases}
\end{equation}
Additionally,
for $x \in \mI_i$ with $i \in \{0,m\}$
\begin{equation}
\label{eqn:NPEr-general-defn2}
N(x,r,c)=
\begin{cases}
\left( Y_{(1)}-r\left( Y_{(1)}-x \right), Y_{(1)} \right)     & \text{if $x < Y_{(1)}$,}\\
\left( Y_{(m)}, Y_{(m)}+r\,\left( x-Y_{(m)} \right) \right) & \text{if $x > Y_{(m)}$.}
\end{cases}
\end{equation}
Notice that for $i \in \{0,m\}$,
the proportional-edge proximity region does not depend on the centrality parameter $c$.
For $x \in \Y_m$,
we define $N(x,r,c)=\{x\}$ for all $r \ge 1$
and if $x = M_{c,i}$, then in Equation \eqref{eqn:NPEr-general-defn1},
we arbitrarily assign $N(x,r,c)$ to be one of
$\left( Y_{(i)}, Y_{(i)}+r\,\left( x-Y_{(i)} \right) \right)\cap \mI_i$ or
$\left( Y_{(i+1)}-r\left( Y_{(i+1)}-x \right), Y_{(i+1)} \right)\cap \mI_i$.
For $c = 0$,
we have $\left( M_{c,i},Y_{(i+1)} \right)=\mI_i$
and
for $c = 1$,
we have $(Y_{(i)},M_{c,i})=\mI_i$.
So,
we set
$N(x,r,0):= \left( Y_{(i+1)}-r\left(Y_{(i+1)}-x\right), Y_{(i+1)} \right)  \cap \mI_i$
and
$N(x,r,1):= \left( Y_{(i)}, Y_{(i)}+r\,\left( x-Y_{(i)} \right) \right) \cap \mI_i$.
For $r > 1$, we have $x \in N(x,r,c)$ for all $x \in \mI_i$.
Furthermore,
$\lim_{r \rightarrow \infty} N(x,r,c) = \mI_i$
for all $x \in \mI_i$,
so we define $N(x,\infty,c) = \mI_i$ for all such $x$.

For $X_i \stackrel{iid}{\sim} F$,
with the additional assumption
that the non-degenerate one-dimensional
probability density function (pdf) $f$ exists with support $\mS(F) \subseteq \mI_i$
and $f$ is continuous around $M_{c,i}$ and around the end points of $\mI_i$,
implies that the special cases in the construction
of $N(\cdot,r,c)$ ---
$X$ falls at $M_{c,i}$ or the end points of $\mI_i$ ---
occurs with probability zero.
For such an $F$,
the region $N(X_i,r,c)$ is an interval a.s.

The data-random proportional-edge PCD has the vertex set $\X_n$ and arc set $\A$ defined by
$(X_i,X_j) \in \A$ iff $X_j \in N(X_i,r,c)$.
We call such digraphs $\D_{n,m}(r,c)$-digraphs.
A $\D_{n,m}(r,c)$-digraph is a {\em pseudo digraph} according some authors,
if loops are allowed (see, e.g., \cite{chartrand:1996}).
The $\D_{n,m}(r,c)$-digraphs are closely related to the {\em proximity graphs} of
\cite{jaromczyk:1992} and might be considered as a special case of
{\em covering sets} of \cite{tuza:1994} and {\em intersection digraphs} of \cite{sen:1989}.
Our data-random proximity digraph is a {\em vertex-random digraph} and
is not a standard random graph (see, e.g., \cite{janson:2000}).
The randomness of a $\D_{n,m}(r,c)$-digraph lies in the fact that
the vertices are random with the joint distribution $F_{X,Y}$,
but arcs $(X_i,X_j)$ are
deterministic functions of the random variable $X_j$ and the random set $N(X_i,r,c)$.
In $\R$, the data-random PCD is a special case of
{\em interval catch digraphs} (see, e.g., \cite{sen:1989} and \cite{prisner:1994}).
Furthermore, when $r=2$ and $c=1/2$ (i.e., $M_{c,i}=\left( Y_{(i)}+Y_{(i+1)} \right)/2$)
we have $N(x,r,c)=B(x,r(x))$
where $B(x,r(x))$ is the ball centered at $x$ with radius $r(x)=d(x,\Y_m)=\min_{y \in \Y_m}d(x,y)$.
The region $N(x,2,1/2)$ corresponds to the proximity region
which gives rise to the CCCD of \cite{priebe:2001}.


\subsection{Domination Number of Random $\D_{n,m}(r,c)$-digraphs}
In a digraph $D=(\V,\A)$ of order $|\V|=n$, a vertex $v$ {\em dominates}
itself and all vertices of the form $\{u:\,(v,u) \in \A\}$.
A {\em dominating set}, $S_D$, for the digraph $D$ is a subset of
$\V$ such that each vertex $v \in \V$ is dominated by a vertex in $S_D$.
A {\em minimum dominating set},  $S^*_D$, is a dominating set of minimum cardinality;
and the {\em domination number}, denoted $\g(D)$, is defined as $\g(D):=|S^*_D|$,
where $|\cdot|$ is the set cardinality functional (\cite{west:2001}).
If a minimum dominating set consists of only one vertex,
we call that vertex a {\em dominating vertex}.
The vertex set $\V$ itself is always a dominating set,
so $\g(D) \le n$.

Let $\F\left( \R^d \right):=\{F_{X,Y} \text{ on } \R^d \text { with } P(X=Y)=0\}$.
As in \cite{priebe:2001} and \cite{ceyhan:dom-num-CCCD-NonUnif},
we consider $\D_{n,m}(r,c)$-digraphs for which
$\X_n$ and $\Y_m$ are random samples from $F_X$ and $F_Y$, respectively,
and the joint distribution of $X,Y$ is $F_{X,Y} \in \F\left( \R^d \right)$.
We call such digraphs \emph{$\F\left( \R^d \right)$-random $\D_{n,m}(r,c)$-digraphs}
and focus on the random variable $\g(D)$.
To make the dependence on sample sizes $n$ and $m$,
the distribution $F$,
and  the parameters $r$ and $c$ explicit,
we use $\g_{{}_{n,m}}(F,r,c)$ instead of $\g(D)$.
For $n \ge 1$ and $m \ge 1$,
it is trivial to see that $1 \le \g_{{}_{n,m}}(F,r,c) \le n$,
and $1 \le \g_{{}_{n,m}}(F,r,c) < n$ for nontrivial digraphs.


\subsection{Special Cases for the Distribution of the
Domination Number of $\F(\R)$-random $\D_{n,m}(r,c)$-digraphs}
\label{sec:domination-number-Dnm}
Let $\X_n$ and $\Y_m$ be two samples from $\F(\R)$, $\X_{[i]}:=\X_n \cap \mI_i$,
and $\Y_{[i]}:=\{Y_{(i)},Y_{(i+1)}\}$ for $i=0,1,2,\ldots,m$.
This yields a disconnected digraph with subdigraphs
each of which might be null or itself disconnected.
Let $D_{[i]}$ be the component of the random $\D_{n,m}(r,c)$-digraph
induced by the pair $\X_{[i]}$ and $\Y_{[i]}$ for $i=0,1,2,\ldots,m$,
$n_i:=\left|\X_{[i]}\right|$, and $F_i$ be the density $F_X$ restricted to $\mI_i$,
and
$\g_{{}_{n_i,2}}(F_i,r,c)$ be the domination number of $D_{[i]}$.
Let also that $M_{c,i} \in \mI_i$ be the point that divides
the interval $\mI_i$ in ratios $c$ and $1-c$
(i.e., length of the subinterval to the left of $M_{c,i}$ is
$c \times 100$ \% of the length of $\mI_i$).
Then $\g_{{}_{n,m}}(F,r,c)=\sum_{i=0}^m ( \g_{{}_{n_i,2}}(F_i,r,c) \I(n_i>0) )$
where $\I(\cdot)$ is the indicator function.
We study the simpler random variable $\g_{{}_{n_i,2}}(F_i,r,c)$ first.
The following lemma follows trivially.

\begin{lemma}
\label{lem:end-intervals}
For $i \in \{ 0,m \}$,
we have $\g_{{}_{n_i,2}}(F_i,r,c)=\I(n_i >0)$ for all $r \ge 1$.
\end{lemma}

Let $\G_1\left( B,r,c \right)$ be the $\G_1$-region
for set $B$ associated with the proximity map on $N(\cdot,r,c)$.

\begin{lemma}
\label{lem:G1-region-in-Ii}
The $\G_1$-region for $\X_{[i]}$ in $\mI_i$ with $r \ge 1$ and $c \in [0,1]$
is
$$\G_1\left( \X_{[i]},r,c \right) =
\Biggl(Y_{(i)}+\frac{\max\,\left( \X_{[i]} \right)}{r},M_{c,i} \Biggr] \bigcup
\Biggl[M_{c,i}, Y_{(i+1)}-\frac{Y_{(i+1)}-\min\left( \X_{[i]} \right)}{r}\Biggr)$$
with the understanding that the intervals $(a,b)$, $(a,b]$, and $[a,b)$ are empty if $a \ge b$.
\end{lemma}
\noindent {\bf Proof:}
By definition,
$\G_1\left( \X_{[i]},r,c \right) = \{x \in \mI_i: \X_{[i]} \subset N(x,r,c)\}$.
Suppose $r \ge 1$ and $c \in [0,1]$.
Then
for $x \in ( Y_{(i)},M_{c,i} ]$,
we have
$\X_{[i]} \subset N(x,r,c)$
iff $Y_{(i)}+r\,(x-Y_{(i)}) > \max\,\left( \X_{[i]} \right)$
iff $x > Y_{(i)}+\frac{\max\,\left( \X_{[i]} \right)}{r}$.
Likewise
for $x \in [ M_{c,i},Y_{(i+1)} )$,
we have
$\X_{[i]} \subset N(x,r,c)$
iff $Y_{(i+1)}-r\,(Y_{(i+1)}-x) < \min\,\left( \X_{[i]} \right)$
iff $x < Y_{(i+1)}-\frac{Y_{(i+1)}-\min\left( \X_{[i]} \right)}{r}$.
Therefore
$\G_1\left( \X_{[i]},r,c \right) =
\Bigl(Y_{(i)}+\frac{\max\,\left( \X_{[i]} \right)}{r},M_{c,i} \Bigr] \bigcup
\Bigl[M_{c,i}, Y_{(i+1)}-\frac{Y_{(i+1)}-\min\left( \X_{[i]} \right)}{r}\Bigr).$
$\blacksquare$

Notice that if $\X_{[i]} \cap \G_1\left( \X_{[i]},r,c \right) \not=\emptyset$,
we have $\g_{{}_{n_i,2}}(F_i,r,c)=1$,
hence the name $\G_1$-region and the notation $\G_1(\cdot)$.
For $i=1,2,3,\ldots,(m-1)$ and $n_i > 0$,
we prove that $\g_{{}_{n_i,2}}(F_i,r,c) = 1$ or $2$
with distribution dependent probabilities.
Hence,
to find the distribution of $\g_{{}_{n_i,2}}(F_i,r,c)$,
it suffices to find $P(\g_{{}_{n_i,2}}(F_i,r,c)=1)$ or $p_{{}_{n_i}}(F_i,r,c):=P\bigl( \g_{{}_{n_i,2}}(F_i,r,c)=2 \bigr)$.
For computational convenience, we employ the latter in our calculations, henceforth.


\begin{theorem}
\label{thm:gamma 1 or 2}
For $i=1,2,3,\ldots,(m-1)$,
let the support of $F_i$ have a positive Lebesgue measure.
Then for $n_i > 1$, $r \in (1,\infty)$, and $c \in (0,1)$,
we have $\g_{{}_{n_i,2}}(F_i,r,c) \sim 1+\BER\left( p_{{}_{n_i}}(F_i,r,c) \right)$.
Furthermore,
$\g_{{}_{1,2}}(F_i,r,c)=1$ for all $r \ge 1$ and $c \in [0,1]$;
$\g_{{}_{n_i,2}}(F_i,r,0)=\g_{{}_{n_i,2}}(F_i,r,1)=1$ for all $n_i \ge 1$ and $ r \ge 1$;
and
$\g_{{}_{n_i,2}}(F_i,\infty,c)=1$ for all $n_i \ge 1$ and $c \in [0,1]$.
\end{theorem}

\noindent {\bf Proof:}
Let $X^-_i:=\argmin_{x \in \X_{[i]} \cap \left( Y_{(i)},M_{c,i} \right)}d(x,M_{c,i})$
provided that $\X_{[i]} \cap \left( Y_{(i)},M_{c,i} \right) \not= \emptyset$,
and $X^+_i:=\argmin_{x \in \X_{[i]} \cap \left( M_{c,i},Y_{(i+1)} \right)}d(x,M_{c,i})$
provided that $\X_{[i]} \cap \left( M_{c,i},Y_{(i+1)} \right) \not= \emptyset$.
That is, $X^-_i$ and $X^+_i$ are closest class $\X$ points (if they exist)
to $M_{c,i}$ from left and right, respectively.
Notice that since $n_i > 0$, at least one of $X^-_i$ and $X^+_i$ exists a.s.
If $\X_{[i]} \cap \left( Y_{(i)},M_{c,i} \right) = \emptyset$,
then $\X_{[i]} \subset N\left( X^+_i,r,c \right)$;
so $\g_{{}_{n_i,2}}(F_i,r,c)=1$.
Similarly,
if $\X_{[i]} \cap (M_{c,i},Y_{(i)})= \emptyset$,
then $\X_{[i]} \subset N\left( X^-_i,r,c \right)$;
so $\g_{{}_{n_i,2}}(F_i,r,c)=1$.
If both of $\X_{[i]} \cap \left( Y_{(i)},M_{c,i} \right)$ and $\X_{[i]} \cap (M_{c,i},Y_{(i)})$ are nonempty,
then $\X_{[i]} \subset N\left( X^-_i,r,c \right) \cup N\left( X^+_i,r,c \right)$,
so $\g_{{}_{n_i,2}}(F_i,r,c) \le 2$.
Since $n_i > 0$, we have $1 \le \g_{{}_{n_i,2}}(F_i,r,c) \le 2$.
The desired result follows,
since the probabilities $1-p_{{}_{n_i}}(F,r,c))=P(\g_{{}_{n_i,2}}(F_i,r,c)=1)$ and
$p_{{}_{n_i}}(F,r,c))=P(\g_{{}_{n_i,2}}(F_i,r,c)=2)$
are both positive.
The special cases in the theorem follow by construction.
$\blacksquare$

The probability $p_{{}_{n_i}}(F,r,c))=P\left( \X_{[i]} \cap \G_1\left( \X_{[i]},r,c \right) =\emptyset \right)$
depends on the conditional distribution $F_{X|Y}$ and the interval $\G_1\left( \X_{[i]},r,c \right)$,
which, if known, will make possible the calculation of $p_{{}_{n_i}}(F_i,r,c)$.
As an immediate result of Lemma \ref{lem:end-intervals} and Theorem \ref{thm:gamma 1 or 2},
we have the following upper bound for $\g_{{}_{n,m}}(F,r,c)$.

\begin{theorem}
\label{thm:gamma-Dnm-r-M}
Let $D_{n,m}(r,c)$ be an $\F(\R)$-random $\D_{n,m}(r,c)$-digraph
and $k_1$, $k_2$, and $k_3$ be three natural numbers defined as
$k_1:=\sum_{i=1}^{m-1} \I(n_i>1)$,
$k_2:=\sum_{i=1}^{m-1} \I(n_i=1)$,
and
$k_3:=\sum_{i \in \{0,m\}} \I(n_i > 0)$.
Then for $n \ge 1,\,m \ge 1$, $r \ge 1$, and $c \in [0,1]$,
we have
$1 \le \g_{{}_{n,m}}(F,r,c) \le 2\,k_1+k_2+k_3 \le \min(n,2\,m)$.
Furthermore,
$\g_{{}_{1,m}}(F,r,c)=1$ for all $m \ge 1$, $r \ge 1$, and $c \in [0,1]$;
$\g_{{}_{n,1}}(F,r,c)=\sum_{i \in \{0,1\}} \I(n_i > 0)$ for all $n \ge 1$ and $r \ge 1$;
$\g_{{}_{1,1}}(F,r,c)=1$ for all $r \ge 1$;
$\g_{{}_{n,m}}(F,r,0)=\g_{{}_{n,m}}(F,r,1)=k_1+k_2+k_3$ for all $m > 1$, $n \ge 1$, and $r \ge 1$;
and
$\g_{{}_{n,m}}(F,\infty,c)=k_1+k_2+k_3$ for all $m > 1$, $n \ge 1$, and $c \in [0,1]$.
\end{theorem}

\noindent {\bf Proof:}
Suppose $n \ge 1,\,m \ge 1$, $r \ge 1$, and $c \in [0,1]$.
Then for $i = 1,2,\ldots,(m-1)$,
by Theorem \ref{thm:gamma 1 or 2},
we have $\g_{{}_{n_i,2}}(F_i,r,c) \in \{1,2\}$ provided that $n_i>1$,
and $\g_{{}_{1,2}}(F_i,r,c) = 1$.
For $i\in \{0,m\}$,
by Lemma \ref{lem:end-intervals},
we have $\g_{{}_{n_i,2}}(F_i,r,c)= \I(n_i > 0)$.
Since $\g_{{}_{n,m}}(F,r,c)=\sum_{i=0}^m (\g_{{}_{n_i,2}}(F_i,r,c)\I(n_i>0))$,
the desired result follows.
The special cases in the theorem follow by construction.
$\blacksquare$

For $r=1$, the distribution of $\g_{{}_{n_i,2}}(F_i,r,c)$ is simpler and
the distribution of $\g_{{}_{n,m}}(F_i,r,c)$ has simpler upper bounds.




\begin{theorem}
\label{thm:gamma-Dnm-r=1-M}
Let $D_{n,m}(1,c)$ be an $\F(\R)$-random $\D_{n,m}(1,c)$-digraph,
$k_3$ be defined as in Theorem \ref{thm:gamma-Dnm-r-M},
and $k_4$ be a natural number defined as
$k_4:=\sum_{i=1}^{m-1} \left[\I\left( \left|\X_{[i]} \cap \left( Y_{(i)},M_{c,i} \right)\right|>0 \right)+
\I\left( \left|\X_{[i]} \cap \left( M_{c,i},Y_{(i+1)} \right)\right|>0 \right)\right]$.
Then for $n \ge 1,\,m > 1$, and $c \in [0,1]$,
we have
$1\le \g_{{}_{n,m}}(F,1,c) = k_3+k_4 \le \min(n,2\,m)$.
\end{theorem}

\noindent {\bf Proof:}
Suppose $n \ge 1,\,m > 1$, and $c \in [0,1]$
and let $X^-_i$ and $X^+_i$ be defined as in the proof of Theorem \ref{thm:gamma 1 or 2}.
Then by construction, $\X_{[i]} \cap \left( Y_{(i)},M_{c,i} \right) \subset N\left( X^-_i,1,c \right)$,
but $N\left( X^-_i,1,c \right) \subseteq \left( Y_{(i)},M_{c,i} \right)$.
So $\left[ \X_{[i]} \cap \left( M_{c,i},Y_{(i+1)} \right) \right] \cap N\left( X^-_i,1,c \right) = \emptyset$.
Similarly
$\X_{[i]} \cap \left( M_{c,i},Y_{(i+1)} \right) \subset N\left( X^+_i,1,c \right)$
and $\left[ \X_{[i]} \cap \left( Y_{(i)},M_{c,i} \right) \right] \cap N\left( X^+_i,1,c \right) = \emptyset$.
Then $\g_{{}_{n_i,2}}(F_i,1,c)=1$,
if $\X_{[i]} \subset \left( Y_{(i)},M_{c,i} \right)$
or
$\X_{[i]} \subset \left( M_{c,i},Y_{(i+1)} \right)$,
and $\g_{{}_{n,m}}(F,1,c)=2$,
if $\X_{[i]} \cap \left( Y_{(i)},M_{c,i} \right) \not= \emptyset$
and $\X_{[i]} \cap \left( M_{c,i},Y_{(i+1)} \right) \not= \emptyset$.
Hence for $i=1,2,3,\ldots,(m-1)$,
we have
$\g_{{}_{n,m}}(F,1,c) =
\I\left( \left|\X_{[i]} \cap \left( Y_{(i)},M_{c,i} \right)\right|>0 \right)+
\I\left( \left|\X_{[i]} \cap \left( M_{c,i},Y_{(i+1)} \right)\right|>0 \right)$,
and for $i\in \{0,m\}$,
we have $\g_{{}_{n_i,2}}(F_i,1,c)= \I(n_i > 0)$.
Since $\g_{{}_{n,m}}(F,1,c)=\sum_{i=0}^m (\g_{{}_{n_i,2}}(F_i,1,c)\I(n_i>0))$,
the desired result follows.
$\blacksquare$

Based on Theorem \ref{thm:gamma-Dnm-r=1-M},
we have
$P(\g_{{}_{n_i,2}}(F,1,c)=1)=
P(\X_{[i]} \subset \left( Y_{(i)},M_{c,i} \right))+
P(\X_{[i]} \subset \left( M_{c,i},Y_{(i+1)} \right))$
and
$P(\g_{{}_{n_i,2}}(F,1,c)=2)=
P(\X_{[i]} \cap \left( Y_{(i)},M_{c,i} \right) \not= \emptyset,
\X_{[i]} \cap \left( M_{c,i},Y_{(i+1)} \right) \not= \emptyset)$.

\section{The Distribution of the Domination Number of Proportional-Edge PCDs
for Uniform Data in One Interval}
\label{sec:gamma-dist-uniform}

In the special case of fixed $\Y_2=\{\y_1,\y_2\}$ with $-\infty<\y_1<\y_2<\infty$
and $\X_n =\{X_1,X_2,\ldots,X_n\}$ a random sample from $\U(\y_1,\y_2)$, the uniform distribution on $(\y_1,\y_2)$,
we have a $\D_{n,2}(r,c)$-digraph
for which $F_X=\U(\y_1,\y_2)$. 
We call such digraphs as \emph{$\U(\y_1,\y_2)$-random $\D_{n,2}(r,c)$-digraphs}
and provide the exact distributions of their domination
number for the whole range of $r$ and $c$.
Let $\g_{{}_{n,2}}(\U,r,c)$ 
be the domination number of the PCD based on $N(\cdot,r,c)$ and $\X_n$
and $p_n(\U,r,c):=P(\g_{{}_{n,2}}(\U,r,c)=2)$,
and  $p(\U,r,c):=\lim_{n\rightarrow \infty}p_n(\U,r,c)$.
We present a ``scale invariance" result for $N(\cdot,r,c)$.
This invariance property will simplify the notation and calculations in
our subsequent analysis by allowing us to consider the special case
of the unit interval, $(0,1)$.

\begin{proposition}
\label{prop:scale-inv-NYr}
(Scale Invariance Property)
Suppose $\X_n$ is a random sample (i.e., a set of iid random variables) from $\U(\y_1,\y_2)$.
Then for any $r \in [1,\infty]$ the distribution of $\g_{{}_{n,2}}(\U,r,c)$ is
independent of $\Y_2$ and hence the support interval $(\y_1,\y_2)$.
\end{proposition}

\noindent {\bf Proof:}
Let $\X_n$ be a random sample from $\U(\y_1,\y_2)$ distribution.
Any $\U(\y_1,\y_2)$ random variable can be transformed into a $\U(0,1)$
random variable by the transformation $\phi(x)=(x-\y_1)/(\y_2-\y_1)$,
which maps intervals $(t_1,t_2) \subseteq (\y_1,\y_2)$ to
intervals $\bigl( \phi(t_1),\phi(t_2) \bigr) \subseteq (0,1)$.
That is,
if $X \sim \U(\y_1,\y_2)$,
then we have $\phi(X) \sim \U(0,1)$
and
$P(X \in (t_1,t_2))=P(\phi(X) \in \bigl( \phi(t_1),\phi(t_2) \bigr)$
for all $(t_1,t_2) \subseteq (\y_1,\y_2)$.
So, without loss of generality, we can assume $\X_n$ 
is a random sample from the $\U(0,1)$ distribution.
Therefore, the distribution of $\g_{{}_{n,2}}(\U,r,c)$
does not depend on the support interval $(\y_1,\y_2)$.
$\blacksquare$


Note that scale invariance of $\g_{{}_{n,2}}(F,\infty,c)$ follows trivially
for all $\X_n$ from any $F$ with support in $(\y_1,\y_2)$,
since for $r=\infty$, we have $\g_{{}_{n,2}}(F,\infty,c)=1$  a.s.
for all $n>1$ and $c \in (0,1)$.
The scale invariance of $\g_{{}_{1,2}}(F,r,c)$ holds
for $n=1$ for all $r \ge 1$ and $c \in [0,1]$,
and scale invariance of $\g_{{}_{n,2}}(F,r,c)$ with $c \in \{0,1\}$ holds
for $n \ge 1$ and $r \ge 1$,  as well.
Based on Proposition \ref{prop:scale-inv-NYr},
for uniform data,
we may assume that $(\y_1,\y_2)$ is the
unit interval $(0,1)$ for $N(\cdot,r,c)$ with general $c$.
Then the proportional-edge proximity region for $x \in (0,1)$
with parameters $r \ge 1$ and $c \in [0,1]$ becomes
\begin{equation}
\label{eqn:NPEr-(0,1)-defn1}
N(x,r,c)=
\begin{cases}
(0, r\,x) \cap (0,1) & \text{if $x \in (0,c)$,}\\
(1-r(1-x), 1)  \cap (0,1)     & \text{if $x \in (c,1)$.}
\end{cases}
\end{equation}

The region $N(c,r,c)$ is arbitrarily taken to be one of
$(0, r\,x) \cap (0,1)$ or $(1-r(1-x), 1)\cap (0,1)$.
Moreover,
$N(0,r,c):=\{0\}$ and $N(1,r,c):=\{1\}$ for all $r \ge 1$ and $c \in [0,1]$;
For $X_i \stackrel{iid}{\sim} \U(\y_1,\y_2)$,
the special cases in the construction of $N(\cdot,r,c)$ ---
$X$ falls at $c$ or the end points of $(\y_1,\y_2)$ ---
occur with probability zero.
Moreover, the region $N(x,r,c)$ is an interval a.s.

The $\G_1$-region, $\G_1(\X_n,r,c)$, depends on $X_{(1)}$, $X_{(n)}$, $r$, and $c$.
If $\G_1(\X_n,r,c) \not= \emptyset$,
then we have $\G_1(\X_n,r,c)=(\delta_1,\delta_2)$
where at least one end points $\delta_1,\delta_2$ is a function of $X_{(1)}$ and $X_{(n)}$.
For $\U(0,1)$ data,
given $X_{(1)}=x_1$ and $X_{(n)}=x_n$,
the probability of $p_n(\U,r,c)$ is
$\displaystyle \left( 1-(\delta_2-\delta_1)/(x_n-x_1)\right)^{(n-2)}$
provided that $\G_1(\X_n,r,c) \not= \emptyset$;
and if $\G_1(\X_n,r,c) = \emptyset$,
then $\g_{{}_{n,2}}(\U,r,c)=2$ holds.
Then
\begin{equation}
\label{eqn:Pg2-int-first}
P(\g_{{}_{n,2}}(\U,r,c)=2,\;\G_1(\X_n,r,c) \not= \emptyset)=
\int\int_{\mS_1}f_{1n}(x_1,x_n)\left(1-\frac{\delta_2-\delta_1}{x_n-x_1}\right)^{(n-2)}\,dx_ndx_1
\end{equation}
where $\mS_1=\{0<x_1<x_n<1:x_1,x_n \not\in \G_1(\X_n,r,c) \text{ and } \G_1(\X_n,r,c) \not= \emptyset\}$
and $f_{1n}(x_1,x_n)=n(n-1)[x_n-x_1]^{(n-2)}\I(0<x_1<x_n<1)$.
The integral simplifies to
\begin{equation}
\label{eqn:Pg2-integral-G1nonempty-unif}
P(\g_{{}_{n,2}}(\U,r,c)=2,\;\G_1(\X_n,r,c) \not= \emptyset)=
\int\int_{\mS_1}n(n-1)[x_n-x_1+\delta_1-\delta_2]^{(n-2)}\,dx_ndx_1.
\end{equation}
If $\G_1(\X_n,r,c) = \emptyset$,
then $\g_{{}_{n,2}}(\U,r,c)=2$.
So
\begin{equation}
\label{eqn:Pg2-integral-G1empty-unif}
P(\g_{{}_{n,2}}(\U,r,c)=2,\;\G_1(\X_n,r,c)= \emptyset)=
\int\int_{\mS_2}f_{1n}(x_1,x_n)\,dx_ndx_1
\end{equation}
where $\mS_2=\{0<x_1<x_n<1:\G_1(\X_n,r,c) = \emptyset\}$.

The probability $p_n(\U,r,c)$ is the sum of the probabilities in
Equations \eqref{eqn:Pg2-integral-G1nonempty-unif} and \eqref{eqn:Pg2-integral-G1empty-unif}.

\subsection{The Exact Distribution of the Domination Number
of $\U(\y_1,\y_2)$-random $\D_{n,2}(2,1/2)$-digraphs}
\label{sec:r=2-and-M_c=1/2}
For $r=2$ and $c=1/2$,
we have
$N(x,2,1/2)=B(x,r(x))$ where $r(x)=\min(x,1-x)$ for $x \in (0,1)$.
Hence proportional-edge PCD based on $N(x,2,1/2)$
is equivalent to the CCCD of \cite{priebe:2001}.
Moreover, $\G_1(\X_n,2,1/2)=\left(X_{(n)}/2,\left(1+X_{(1)}\right)/2\right)$.
It has been shown that $p_n(\U,2,1/2)=4/9-(16/9) \, 4^{-n}$
(\cite{priebe:2001}).
Hence, for $\U(\y_1,\y_2)$ data with $n \ge 1$, we have
\begin{equation}
\label{eqn:finite-sample-unif}
\g_{{}_{n,2}}(\U,2,1/2)=  \left\lbrace \begin{array}{ll}
       1           & \text{w.p. $5/9+(16/9) \, 4^{-n}, $}\\
       2           & \text{w.p. $4/9-(16/9) \, 4^{-n}, $}
\end{array} \right.
\end{equation}
where w.p. stands for ``with probability".
Then as $n \rightarrow \infty$,
$\g_{{}_{n,2}}(\U,2,1/2)$ converges in distribution to $1+\BER(4/9)$.
For $m>2$, \cite{priebe:2001} computed the exact distribution of $\g_{{}_{n,m}}(\U,2,1/2)$ also.
However, 
the scale invariance property does not hold for general $F$;
that is, for $X_i \stackrel{iid}{\sim}F$ with support $\mS(F) \subseteq(\y_1,\y_2)$,
the exact and asymptotic distribution of $\g_{{}_{n,2}}(F,2,1/2)$
depends on $F$ and $\Y_2$ (\cite{ceyhan:dom-num-CCCD-NonUnif}).

\subsection{The Exact Distribution of the Domination Number of $\U(\y_1,\y_2)$-random $\D_{n,2}(2,c)$-digraphs}
\label{sec:r=2-and-M_c=c}
For $r=2$, $c \in (0,1)$, and $(\y_1,\y_2)=(0,1)$,
the $\G_1$-region is $\G_1(\X_n,2,c)=( X_{(n)}/2,c ] \cup [ c,(1+X_{(1)})/2)$.
Notice that $( X_{(n)}/2,c ]$ or $[ c,(1+X_{(1)})/2 )$
could be empty, but not simultaneously.
\begin{theorem}
\label{thm:r=2 and M_c=c}
For $\U(\y_1,\y_2)$ data
and $n \ge 1$,
we have
$\g_{{}_{n,2}}(\U,2,c) \sim 1+\BER(p_{{}_n}(\U,2,c))$
where
$p_{{}_n}(\U,2,c) =
\nu_{1,n}(c)\I(c \in (0,1/3]+
\nu_{2,n}(c)\I(c \in (1/3,1/2]+
\nu_{3,n}(c)\I(c \in (1/2,2/3]+
\nu_{4,n}(c)\I(c \in (2/3,1)
$
with
$$\nu_{1,n}(c)= \frac{2}{3} \left(c+\frac{1}{2}\right)^n-\frac{8}{9} 4^{-n}-\frac{2}{3}\left(\frac{1-c}{2}\right)^n
+\frac{1}{9}(1-3c)^n-\frac{2}{9}\left(3c-\frac{1}{2}\right)^n,$$
$$\nu_{2,n}(c)= \frac{2}{3} \left(c+\frac{1}{2}\right)^n-\frac{8}{9} 4^{-n}-\frac{2}{3}\left(\frac{1-c}{2}\right)^n-
\frac{2}{9}\left(\frac{3c-1}{2}\right)^n-\frac{2}{9}\left(3c-\frac{1}{2}\right)^n,$$
$\nu_{3,n}(c)=\nu_{2,n}(1-c)$, and $\nu_{4,n}(c)=\nu_{1,n}(1-c)$.
Furthermore,
$\g_{{}_{n,2}}(\U,2,0)=\g_{{}_{n,2}}(\U,2,1)=1$ for all $n \ge 1$.
\end{theorem}

Observe that
the parameter $p_{{}_n}(\U,2,c)$ is continuous in $c \in (0,1)$ for fixed $n < \infty$,
but there are jumps (hence discontinuities) in $p_{{}_n}(\U,2,c)$ at $c \in \{0,1\}$.
In particular,
$\lim_{c \rightarrow 0}p_{{}_n}(\U,2,c)=
\lim_{c \rightarrow 1}p_{{}_n}(\U,2,c)=
\lim_{c \rightarrow 0}\nu_{1,n}(c)=
\lim_{c \rightarrow 1}\nu_{4,n}(c)=
\frac{1}{9}-\frac{2}{9}(-2)^n-\frac{8}{9}4^{-n}$,
but
$p_{{}_n}(\U,2,0)=p_{{}_n}(\U,2,1)=0$ for all $n \ge 1$.
For $c=1/2$,
we have $p_{{}_n}(\U,2,c)=4/9-(16/9) \, 4^{-n}$,
hence the distribution of $\g_{{}_{n,2}}(\U,2,c=1/2)$ is same as in Equation \eqref{eqn:finite-sample-unif}.

In the limit
as $n \rightarrow \infty$,
for $c \in [0,1]$,
we have
\begin{equation*}
\g_{{}_{n,2}}(\U,2,c) \sim
\left\lbrace \begin{array}{ll}
       1+\BER(4/9), & \text{for $c = 1/2$,}\\
       1,           & \text{for $c \not= 1/2$.}\\
\end{array} \right.
\end{equation*}
Observe also the interesting behavior of the asymptotic
distribution of $\g_{{}_{n,2}}(\U,2,c)$ around $c=1/2$.
The parameter $p(\U,2,c)$ is continuous in $c \in [0,1] \setminus \{1/2\}$
(in fact it is unity),
but there is a jump (hence discontinuity) in $p(\U,2,c)$ at $c=1/2$,
since $p(\U,2,1/2)=4/9$
and $p(\U,2,c)=0$ for $c \not= 1/2$.
Hence
for $c = 1/2$, the asymptotic distribution is non-degenerate,
and
for $c \not= 1/2$, the asymptotic distribution is degenerate.
That is, for $c=1/2 \pm \varepsilon$ with $\varepsilon>0$ arbitrarily small,
although the exact distribution is non-degenerate,
the asymptotic distribution is degenerate.

\subsection{The Exact Distribution of the Domination Number of $\U(\y_1,\y_2)$-random $\D_{n,2}(r,1/2)$-digraphs}
\label{sec:r-and-M_c=1/2}
For $r \ge 1$, $c=1/2$, and $(\y_1,\y_2)=(0,1)$,
the $\G_1$-region is $\G_1(\X_n,r,1/2)=(X_{(n)}/r,1/2] \cup [1/2,(r-1+X_{(1)})/r)$
where
$(X_{(n)}/r,1/2]$ or $[1/2,(r-1+X_{(1)})/r)$ could be empty, but not simultaneously.
\begin{theorem}
\label{thm:r and M_c=1/2}
For $\U(\y_1,\y_2)$ data
with $n \ge 1$,
we have
$\g_{{}_{n,2}}(\U,r,1/2) \sim 1+\BER(p_{{}_n}(\U,r,1/2))$
where
\begin{eqnarray*}
p_{{}_n}(\U,r,1/2)=
\begin{cases}
\frac{2\,r}{(r+1)^2} \left( \left(\frac{2}{r}\right)^{n-1}-\left(\frac{r-1}{r^2} \right)^{n-1} \right) &\text{for} \quad r \ge 2, \\
1-\frac{1+r^{2n-1}}{(2\,r)^{n-1}(r+1)}+\frac{(r-1)^n}{(r+1)^2} \left( 1- \left(\frac{r-1}{2\,r}\right)^{n-1}\right)
&\text{for} \quad 1 \le r < 2.
\end{cases}
\end{eqnarray*}
\end{theorem}

Notice that for fixed $n < \infty$,
the parameter
$p_{{}_n}(\U,r,1/2)$ is continuous in $r \ge 2$.
In particular,
for $r=2$,
we have $p_{{}_n}(\U,2,1/2)=4/9-(16/9) \, 4^{-n}$,
hence the distribution of $\g_{{}_{n,2}}(\U,r=2,1/2)$ is same as in Equation \eqref{eqn:finite-sample-unif}.
Furthermore,
$\lim_{r \rightarrow 1}p_{{}_n}(\U,r,1/2)=
p_{{}_n}(\U,1,1/2)=
1-2^{1-n}$
and
$\lim_{r \rightarrow \infty}p_{{}_n}(\U,r,1/2)=
p_{{}_n}(\U,\infty,1/2)=
0$.

In the limit,
as $n \rightarrow \infty$, we have
\begin{equation*}
\g_{{}_{n,2}}(\U,r,1/2) \sim
\left\lbrace \begin{array}{ll}
       1           & \text{for $r > 2$,}\\
       1+\BER(4/9) & \text{for $r = 2$,}\\
       2           & \text{for $1 \le r < 2$.}\\
\end{array} \right.
\end{equation*}
Observe the interesting behavior of the asymptotic distribution
of $\g_{{}_{n,2}}(\U,r,1/2)$ around $r=2$.
The parameter
$p(\U,r,1/2)$ is continuous (in fact piecewise constant) for $r \in [1,\infty) \setminus \{2\}$.
Hence for $r \not= 2$, the asymptotic distribution is degenerate,
as $p(\U,r,1/2) = 1$ for $r>2$
and $p(\U,r,1/2) = 2$ w.p. 1 for $r<2$.
That is, for $r=2 \pm \varepsilon$ with $\varepsilon>0$ arbitrarily small,
although the exact distribution is non-degenerate,
the asymptotic distribution is degenerate.

\subsection{The Distribution of the Domination Number of $\U(\y_1,\y_2)$-random $\D_{n,2}(r,c)$-digraphs}
\label{sec:r-and-M}
For $r \ge 1$ and $c \in (0,1)$,
the $\G_1$-region is $\G_1(\X_n,r,c)=(X_{(n)}/r,c] \cup [c,(r-1+X_{(1)})/r)$
where $(X_{(n)}/r,c]$ or $[c,(r-1+X_{(1)})/r)$ could be empty, but not simultaneously.
\begin{theorem}
\label{thm:r and M}
\textbf{Main Result 1:}
For $\U(\y_1,\y_2)$ data with $n \ge 1$, $r \ge 1$, and $c \in ((3-\sqrt{5})/2,1/2)$,
we have
$\g_{{}_{n,2}}(\U,r,c) \sim 1+\BER(p_{{}_n}(\U,r,c))$
where
$p_{{}_n}(\U,r,c) =\pi_{1,n}(r,c) \,\I(r \ge 1/c) + \pi_{2,n}(r,c) \,\I(1/(1-c) \le r < 1/c)+
\pi_{3,n}(r,c) \,\I((1-c)/c \le r < 1/(1-c)) + \pi_{4,n}(r,c) \,\I( 1 \le r < (1-c)/c)
$
with
$$\pi_{1,n}(r,c)=\frac{2\,r}{(r+1)^2} \left( \left(\frac{2}{r}\right)^{n-1}-\left(\frac{r-1}{r^2} \right)^{n-1} \right),$$
\begin{multline*}
\pi_{2,n}(r,c)=
\frac{1}{(r+1)r^{n-1}}\left[ (1+c\,r)^n-(1-c)^n -\frac{1}{r+1}(c\,r^2-r+c\,r+1)^n
-\frac{(r-1)^{n-1}}{r+1}\left(\frac{1}{r^{n-2}}+(c\,r-1+c)^n\right)\right],
\end{multline*}
\begin{multline*}
\pi_{3,n}(r,c)=
1+\frac{(r-1)^{n-1}}{(r+1)^2}
\left[(r-1)-\frac{1}{r^{n-1}}((c\,r-1+c)^n+(r-c\,r-c)^n) \right]
-\frac{1}{r+1}[c^n+(1-c)^n]\left(r^n+\frac{1}{r^{n-1}}\right),
\end{multline*}
and
\begin{multline*}
\pi_{4,n}(r,c)=
1+\frac{(r-1)^{n-1}}{(r+1)^2}(1-c\,r-c)^n\left( r^2-\left(\frac{-1}{r}\right)^{n-1} (r^2-1)\right)+
\frac{(r-1)^n}{(r+1)^2}\left(1-r\left(\frac{r-c\,r-c}{r}\right)^n\right)\\
-\frac{1}{r+1}[c^n+(1-c)^n]\left(r^n-\frac{1}{r^{n-1}}\right).
\end{multline*}
And for $c \in (0,(3-\sqrt{5})/2]$,
we have
$p_{{}_n}(\U,r,c) =\vartheta_{1,n}(r,c) \,\I(r \ge 1/c) + \vartheta_{2,n}(r,c) \,\I((1-c)/c \le r < 1/c)+
\vartheta_{3,n}(r,c) \,\I(1/(1-c) \le r < (1-c)/c) + \vartheta_{4,n}(r,c) \,\I( 1 \le r < 1/(1-c))
$
where
$\vartheta_{1,n}(r,c)=\pi_{1,n}(r,c)$,
$\vartheta_{2,n}(r,c)=\pi_{2,n}(r,c)$,
$\vartheta_{4,n}(r,c)=\pi_{4,n}(r,c)$,
and
\begin{multline*}
\vartheta_{3,n}(r,c)=
\frac{r}{(r+1)^2}\Biggl[
(r-1)^{n-1}(1-c\,r-c)^n\left(r+(r^2-1)\left(\frac{-1}{r}\right)^n\right)-
\left(\frac{r-1}{r^2}\right)^{n-1}-
\left(\frac{c\,r^2-c+c\,r+1}{r}\right)^n+\\
\frac{r+1}{r^n}[(1+c\,r)^n-(1-c)^n]\Biggr].
\end{multline*}
Furthermore, we have
$\g_{{}_{n,2}}(\U,r,0)=\g_{{}_{n,2}}(\U,r,1)=1$ for all $n \ge 1$.
\end{theorem}

\begin{figure}[ht]
\centering
\rotatebox{-90}{ \resizebox{2.5 in}{!}{ \includegraphics{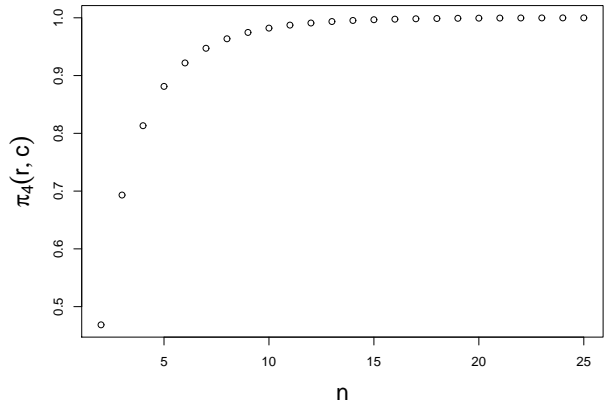}}}
\caption{
\label{fig:pi4-rc}
The probability $\pi_{4,n}(r,c)$ in Main Result 1
with $r=1.2$ and $c=0.4$ for $n=2,3,\ldots,25$.
}
\end{figure}

\begin{figure}[ht]
\centering
\psfrag{teta}{\Huge{$\vartheta_{3,n}(r,c)$}}
\rotatebox{-90}{ \resizebox{2.5 in}{!}{ \includegraphics{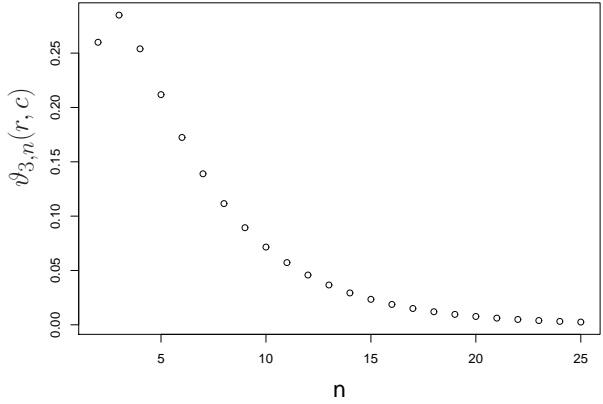}}}
\caption{
\label{fig:teta3-rc}
The probability $\vartheta_{3,n}(r,c)$ in Main Result 1
with $r=2$ and $c=0.3$ for $n=2,3,\ldots,25$.
}
\end{figure}

Some remarks are in order for the Main Result 1.
The partitioning of $c \in (0,1/2)$
as $c \in (0,(3-\sqrt{5})/2]$ and $c \in ((3-\sqrt{5})/2,1/2)$
is due to the relative positions of $1/(1-c)$ and $(1-c)/c$.
For $c \in ((3-\sqrt{5})/2,1/2)$, we have $1/(1-c) > (1-c)/c$
and
for $c \in (0,(3-\sqrt{5})/2)$,
we have $1/(1-c) < (1-c)/c$.
At $c=(3-\sqrt{5})/2$,
$1/(1-c) = (1-c)/c = (\sqrt{5}+1)/2$
and only
$\pi_{1,n}(r,(3-\sqrt{5})/2)=\vartheta_{1,n}(r,(3-\sqrt{5})/2)$,
$\pi_{2,n}(r,(3-\sqrt{5})/2)=\vartheta_{2,n}(r,(3-\sqrt{5})/2)$,
and
$\pi_{4,n}(r,(3-\sqrt{5})/2)=\vartheta_{4,n}(r,(3-\sqrt{5})/2)$ terms survive.
Also, notice the $(-1)^{n}$ terms in $\pi_{4,n}(r,c)$ and $\vartheta_{3,n}(r,c)$
which might suggest fluctuations of these probabilities as $n$ changes (increases).
However,
as $n$ increases,
$\pi_{4,n}(r,c)$ strictly increases towards 1
(see Figure \ref{fig:pi4-rc}),
and
$\vartheta_{3,n}(r,c)$
decreases (strictly decreases for $n \ge 3$) towards 0
(see Figure \ref{fig:teta3-rc}).

\begin{remark}
By symmetry,
in Theorem \ref{thm:r and M},
for $c \in (1/2,(\sqrt{5}-1)/2)$,
we have
$p_{{}_n}(\U,r,c) =\pi_{1,n}(r,1-c) \,\I(r \ge 1/(1-c)) + \pi_{2,n}(r,1-c) \,\I(1/c \le r < 1/(1-c))+
\pi_{3,n}(r,1-c) \,\I(c/(1-c) \le r < 1/c) + \pi_{4,n}(r,1-c) \,\I( 1 \le r < c/(1-c))
$
and
for $c \in [(\sqrt{5}-1)/2),1)$,
$p_{{}_n}(\U,r,c) =\vartheta_{1,n}(r,1-c) \,\I(r \ge 1/(1-c)) + \vartheta_{2,n}(r,1-c) \,\I(c/(1-c) \le r < 1/(1-c))+
\vartheta_{3,n}(r,1-c) \,\I(1/c \le r < c/(1-c)) + \vartheta_{4,n}(r,1-c) \,\I( 1 \le r < 1/c)
$.
$\square$
\end{remark}

Observe that $\lim_{r \rightarrow 1} p_{{}_n}(\U,r,c)=\lim_{r \rightarrow 1} \pi_{4,n}(r,c)=1$
as expected.
For fixed $1 < n < \infty$,
the probability $p_{{}_n}(\U,r,c)$ 
is continuous in $(r,c) \in \{(r,c) \in \R^2: r \ge 1, 0 < c < 1\}$.
In particular,
for $c \in ((3-\sqrt{5})/2,1/2)$,
as $(r,c) \rightarrow (2,1/2)$ in $\{(r,c) \in \R^2: r \ge 1/c\}$,
$p_{{}_n}(\U,r,c)=\pi_{1,n}(r,c) \rightarrow 4/9-(16/9) \, 4^{-n}$;
as $(r,c) \rightarrow (2,1/2)$ in $\{(r,c) \in \R^2: 1/(1-c) \le r < 1/c\}$,
$p_{{}_n}(\U,r,c)=\pi_{2,n}(r,c) \rightarrow 4/9-(16/9) \, 4^{-n}$;
and
as $(r,c) \rightarrow (2,1/2)$ in $\{(r,c) \in \R^2: (1-c)/c \le r < 1/(1-c)\}$,
$p_{{}_n}(\U,r,c)=\pi_{3,n}(r,c) \rightarrow 4/9-(16/9) \, 4^{-n}$.
The limit $(r,c) \rightarrow (2,1/2)$ is not possible for $\{(r,c) \in \R^2: 1 \le r < (1-c)/c\}$.
For $c \in (0,(3-\sqrt{5})/2]$,
$(r,c) \rightarrow (2,1/2)$ can not occur either.
And for $(r,c)=(2,1/2)$,
the distribution of $\g_{{}_{n,2}}(\U,r,c)$ is $1+\BER(p_{{}_n}(\U,2,1/2))$,
where $p_{{}_n}(\U,2,1/2)=4/9-(16/9) \, 4^{-n}$ as in Equation \eqref{eqn:finite-sample-unif}.
Therefore for fixed $1<n<\infty$,
as $(r,c) \rightarrow (2,1/2)$ in $S=\{(r,c) \in \R^2: r \ge 1, 0 < c <1/2\}$,
we have
$p_{{}_n}(\U,r,c) \rightarrow 4/9-(16/9) \, 4^{-n}$.
Hence
as $(r,c) \rightarrow (2,1/2)$ in $S$,
$\g_{{}_{n,2}}(\U,r,c)$ converges in distribution to $\g_{{}_{n,2}}(\U,2,1/2)$.
However,
$p_{{}_n}(\U,r,c)$ has jumps (hence discontinuities) at $c \in \{0,1\}$.
As $c \rightarrow 0^+$ (which implies
we should consider $c \in (0,(3-\sqrt{5})/2]$,
$1/c \rightarrow \infty$,
$(1-c)/c \rightarrow \infty$,
and
$1/(1-c) \rightarrow 1^+$.
Hence
$\lim_{c\rightarrow 0^+}\vartheta_{1,n}(r,c)=\vartheta_{1,n}(\infty,0)=0$,
$\lim_{c\rightarrow 0^+}\vartheta_{2,n}(r,c)=\vartheta_{2,n}(\infty,0)=0$.
Moreover,
$\lim_{c\rightarrow 0^+}\vartheta_{3,n}(r,c)=\frac{r(r-1)^{n-1}}{(r+1)^2}[r+(-r)^{1-n}(1-r)-r^{2-2n}]$;
$\lim_{c\rightarrow 0^+}\vartheta_{4,n}(r,c)=1+\frac{1}{r+1}[(r-1)^n(1-(-r)^{1-n})-r^n-r^{1-n}]$.
But $p_{{}_n}(\U,r,0)=0$ for all $r \ge 1$.
Similar results can be obtained as $c \rightarrow 1^{-}$.
Observe also that $\lim_{r \rightarrow 1} p_{{}_n}(\U,r,c)=\lim_{r \rightarrow 1} \pi_{4,n}(r,c)=1$.

\begin{theorem}
\label{thm:r and M-asy}
\textbf{Main Result 2:}
Let $D_{n,2}(r,c)$ be based on $\U(\y_1,\y_2)$ data with
$c \in (0,1)$ and $\tau=max(c,1-c)$.
Then the domination number $\g_{{}_{n,2}}(\U,r,c)$ of the PCD has the following
asymptotic distribution.
As $n \rightarrow \infty$,
\begin{equation}
\label{eqn:asy-unif-rM}
\g_{{}_{n,2}}(\U,r,c) \sim
\left\lbrace \begin{array}{ll}
       1+\BER(r/(r+1)),           & \text{for $r = 1/\tau$,}\\
       1,           & \text{for $r> 1/\tau$,}\\
       2,           & \text{for $1 \le r < 1/\tau$.}\\
\end{array} \right.
\end{equation}
\end{theorem}

Notice the interesting behavior of the asymptotic distribution
of $\g_{{}_{n,2}}(\U,r,c)$ around $r=1/\tau$ for any given $c \in (0,1)$.
The asymptotic distribution is non-degenerate only for $r = 1/\tau$.
For $r>1/\tau$,
$\lim_{n \rightarrow \infty}\g_{{}_{n,2}}(\U,r,c)) = 1$ w.p. 1,
and
for $1 \le r < 1/\tau $, $\lim_{n \rightarrow \infty}\g_{{}_{n,2}}(\U,r,1/2)) = 2$ w.p. 1.
%
The critical value $r=1/\tau$ corresponds to
$c=(r-1)/r$, if $c \in(0,1/2)$ (i.e., $\tau=1-c$)
and
$c=1/r$, if $c \in(1/2,1)$ (i.e., $\tau=c$)
and these are only possible for $r \in (1,2)$.
That is, for $r=(1/\tau) \pm \varepsilon$ for $\varepsilon$ arbitrarily small,
although the exact distribution is non-degenerate,
the asymptotic distribution is degenerate.
The parameter $p(\U,r,c)$ is continuous in $r$ and $c$ for $(r,c) \in S \setminus \{1/\tau,c\}$
and there is a jump (hence discontinuity) in the probability $p(\U,r,c)$
at $r=1/\tau$,
since $p(\U,1/\tau,c)=1/(1+\tau)=r/(r+1)$.
Therefore, given a centrality parameter $c \in (0,1)$,
we can choose the expansion parameter $r$
for which the asymptotic distribution is non-degenerate,
and vice versa.

There is yet another interesting behavior of the asymptotic distribution
around $(r,c)=(2,1/2)$.
The parameter $p(\U,r,c)$
has jumps at $c=1/r$ and $(r-1)/r$ for $r \in [1,2]$
with $p(\U,r,1/r)=p(\U,r,(r-1)/r)=r/(r+1)$.
That is,
for fixed $(r,c) \in S$,
$\lim_{n \rightarrow \infty}p_{{}_n}(\U,r,(r-1)/r)=
\lim_{n \rightarrow \infty}p_{{}_n}(\U,r,1/r)= r/(r+1)$.
Letting $(r,c) \rightarrow (2,1/2)$
(i.e., $r \rightarrow 2$)
we get
$p(\U,r,(r-1)/r) \rightarrow 2/3$
and
$p(\U,r,1/r) \rightarrow 2/3$,
but
$p(\U,2,1/2)=4/9$.
Hence for $r \in [1,2)$
the distributions
of $\g_{{}_{n,2}}(\U,r,(r-1)/r)$ and $\g_{{}_{n,2}}(\U,r,1/r)$
are identical and both
converge to $1+\BER(r/(r+1))$,
but the distribution of
$\g_{{}_{n,2}}(\U,2,1/2)$
converges to $1+\BER(4/9)$ as $n \rightarrow \infty$.
In other words,
$p(\U,r,(r-1)/r)=p(\U,r,1/r)$
has another jump at $r=2$ (which corresponds to $(r,c)=(2,1/2)$).
This interesting behavior might be due to the symmetry around $c=1/2$.
Because for $c \in (0,1/2)$,
with $r=1/(1-c)$,
for sufficiently large $n$,
a point $X_i$ in $(c,1)$ can dominate all the points in $\X_n$
(implying $\g_{{}_{n,2}}(\U,r,(r-1)/r)=1$),
but no point in $(0,c)$ can dominate all points a.s.
Likewise,
for $c \in (1/2,1)$ with $r=1/c$,
for sufficiently large $n$,
a point $X_i$ in $(0,c)$ can dominate all the points in $\X_n$
(implying $\g_{{}_{n,2}}(\U,r,1/r)=1$),
but no point in $(c,1)$ can dominate all points a.s.
However,
for $c=1/2$ and $r=2$,
for sufficiently large $n$,
points to the left or right of $c$
can dominate all other points in $\X_n$.

\section{The Distribution of the Domination Number for $\F(\R)$-random $\D_{n,2}(r,c)$-digraphs}
\label{sec:non-uniform}
Let $\F(\y_1,\y_2)$ be a family of continuous distributions
with support in $\mS_F \subseteq (\y_1,\y_2)$.
Consider a distribution function $F \in \F(\y_1,\y_2)$.
For simplicity, assume $\y_1=0$ and $\y_2=1$.
Let $\X_n$ be a random sample from $F$,
$\G_1(\X_n,r,c)=(\delta_1,\delta_2)$,
$p_{{}_n}(F,r,c):=P(\g_{{}_{n,2}}(F,r,c)=2)$,
and $p(F,r,c):=\lim_{n \rightarrow \infty}P(\g_{{}_{n,2}}(F,r,c)=2)$.
The exact (i.e., finite sample) and asymptotic distributions of $\g_{{}_{n,2}}(F,r,c)$
are $1+\BER\left(p_{{}_n}(F,r,c)\right)$ and $1+\BER\left(p(F,r,c)\right)$, respectively.
That is, for finite $n > 1$, $r \in [1,\infty)$, and $c \in (0,1)$,
we have
\begin{equation}
\g_{{}_{n,2}}(F,r,c)=  \left\lbrace \begin{array}{ll}
       1           & \text{w.p. $1-p_{{}_n}(F,r,c)$},\\
       2           & \text{w.p. $p_{{}_n}(F,r,c)$}.
\end{array} \right.
\end{equation}
Moreover,
$\g_{{}_{1,2}}(F,r,c)=1$ for all $r \ge 1$ and $c \in [0,1]$,
$\g_{{}_{n,2}}(F,r,0)=\g_{{}_{1,2}}(F,r,1)=1$ for all $n \ge 1$ and $r \ge 1$,
$\g_{{}_{n,2}}(F,\infty,c)=1$ for all $n \ge 1$ and $c \in [0,1]$,
and
$\g_{{}_{n,2}}(F,1,c)=k_4$ for all $n \ge 1$ and $c \in (0,1)$
where $k_4$ is as in Theorem \ref{thm:gamma-Dnm-r=1-M} with $m=2$.
The asymptotic distribution is similar with
$p_{{}_n}(F,r,c)$ being replaced by $p(F,r,c)$.
The special cases are similar in the asymptotics
with the exception that
$p(F,1,c)=1$ for all $c \in (0,1)$.
The finite sample mean and variance of $\g_{{}_{n,2}}(F,r,c)$ are given by
$1+p_{{}_n}(F,r,c)$ and $p_{{}_n}(F,r,c)\,(1-p_{{}_n}(F,r,c))$, respectively;
and the asymptotic mean and variance of $\g_{{}_{n,2}}(F,r,c)$ are given by
$1+p(F,r,c)$ and $p(F,r,c)\,(1-p(F,r,c))$, respectively.

Given $X_{(1)}=x_1$ and $X_{(n)}=x_n$,
the probability of $\g_{{}_{n,2}}(F,r,c)=2$ (i.e., $p_{{}_n}(F,r,c)$) is
$\displaystyle ( 1-[F(\delta_2)-F(\delta_1)]/[F(x_n)-F(x_1)])^{(n-2)}$
provided that $\G_1(\X_n,r,c)=(\delta_1,\delta_2) \not= \emptyset$;
if $\G_1(\X_n,r,c) = \emptyset$,
then $\g_{{}_{n,2}}(F,r,c)=2$.
Then
\begin{equation}
\label{eqn:Pg2-int-first}
P(\g_{{}_{n,2}}(F,r,c)=2,\;\G_1(\X_n,r,c) \not= \emptyset)=
\int\int_{\mS_1}f_{1n}(x_1,x_n)\left(1-\frac{F(\delta_2)-F(\delta_1)}{F(x_n)-F(x_1)}\right)^{(n-2)}\,dx_ndx_1
\end{equation}
where $\mS_1=\{0<x_1<x_n<1:(x_1,x_n) \not\in \G_1(\X_n,r,c),\G_1(\X_n,r,c) \not= \emptyset\}$
and $f_{1n}(x_1,x_n)=n(n-1)f(x_1)f(x_n)\bigl[F(x_n)-F(x_1)\bigr]^{(n-2)}\I(0<x_1<x_n<1)$
which is the joint pdf of $X_{(1)},X_{(n)}$.
The integral in \eqref{eqn:Pg2-int-first} becomes
\begin{equation}
\label{eqn:Pg2-integral-G1nonempty-U}
P(\g_{{}_{n,2}}(F,r,c)=2,\;\G_1(\X_n,r,c) \not= \emptyset)=\\
\int\int_{\mS_1}H(x_1,x_n)\,dx_ndx_1,
\end{equation}
where
\begin{equation}
\label{eqn:integrand}
H(x_1,x_n):=n\,(n-1)f(x_1)f(x_n)\bigl[F(x_n)+F\left(\delta_1 \right)-
\left( F\left( \delta_2 \right)+F(x_1) \right)\bigr]^{n-2}.
\end{equation}

If $\G_1(\X_n,r,c) = \emptyset$,
then $\g_{{}_{n,2}}(F,r,c)=2$.
So
\begin{equation}
\label{eqn:Pg2-integral-G1empty}
P(\g_{{}_{n,2}}(F,r,c)=2,\;\G_1(\X_n,r,c)= \emptyset)=
\int\int_{\mS_2}f_{1n}(x_1,x_n)\,dx_ndx_1
\end{equation}
where $\mS_2=\{0<x_1<x_n<1:\G_1(\X_n,r,c) = \emptyset\}$.

The probability $p_{{}_n}(F,r,c)$ is the sum of the probabilities
in Equations \eqref{eqn:Pg2-integral-G1nonempty-U} and \eqref{eqn:Pg2-integral-G1empty}.

For $\Y_2=\{\y_1,\y_2\} \subset \R$ with $-\infty<\y_1<\y_2<\infty$,
a quick investigation shows that the $\G_1$-region is
$\G_1(\X_n,r,c)= (\y_1+(X_{(n)}-\y_1)/r,M_c] \cup [M_c,\y_2-(\y_2-X_{(1)})/r)$.
Notice that for a given $c \in [0,1]$,
the corresponding $M_c \in [\y_1,\y_2]$
is $M_c=\y_1+c(\y_2-\y_1)$.


Let $F$ be a continuous distribution with support $\mS(F)\subseteq (0,1)$.
The simplest of such distributions is $\U(0,1)$,
which yields the simplest exact distribution for $\g_{{}_{n,2}}(F,r,c)$
with $(r,c)=(1,1/2)$.
If $X \sim F$, then by probability integral transform, $F(X) \sim \U(0,1)$.
So for any continuous $F$,
we can construct a proximity map depending on $F$ for which
the distribution of the domination number of the
associated digraph has the same distribution
as that of $\g_{{}_{n,2}}(\U,r,c)$.

\begin{proposition}
\label{prop:NF vs NPE}
Let $X_i \stackrel{iid}{\sim} F$ which is an absolutely continuous
distribution with support $\mS(F)=(0,1)$ and $\X_n=\{X_1,X_2,\ldots,X_n\}$.
Define the proximity map $N_F(x,r,c):=F^{-1}(N(F(x),r,c))$.
Then the domination number of the digraph based on $N_F$, $\X_n$, and $\Y_2=\{0,1\}$
has the same distribution as $\g_{{}_{n,2}}(\U,r,c)$.
\end{proposition}
\noindent
{\bfseries Proof:}
Let $U_i:=F(X_i)$ for $i=1,2,\ldots,n$ and $\U_n:=\{U_1,U_2,\ldots,U_n\}$.
Hence, by probability integral transform, $U_i \stackrel{iid}{\sim} \U(0,1)$.
Let $U_{(i)}$ be the $i^{th}$ order statistic of $\U_n$ for $i=1,2,\ldots,n$.
Furthermore,
an absolutely continuous $F$ preserves order;
that is,
for $x \le y$,
we have $F(x) \le F(y)$.
So the image of $N_F(x,r,c)$ under $F$ is
$F(N_F(x,r,c))=N(F(x),r,c)$ for (almost) all $x \in (0,1)$.
Then $F(N_F(X_i,r,c))=N(F(X_i),r,c)=N(U_i,r,c)$ for $i =1,2,\ldots,n$.
Since $U_i \stackrel{iid}{\sim} \U(0,1)$, the distribution of
the domination number of the digraph based on $N(\cdot,r,c)$, $\U_n$, and $\{0,1\}$
is given in Theorem \ref{thm:r and M}.
Observe that for any $j$,
$X_j \in N_F(X_i,r,c)$ iff
$X_j \in F^{-1}(N(F(X_i),r,c))$ iff
$F(X_j) \in N(F(X_i),r,c)$ iff
$U_j \in N(U_i,r,c)$ for $i=1,2,\ldots,n$.
Hence $P(\X_n \subset N_F(X_i,r,c))=P(\U_n \subset N(U_i,r,c)$
for all $i=1,2,\ldots,n$.
Therefore, $\X_n \cap \G_1(\X_n,N_F(r,c))=\emptyset$
iff $\U_n \cap \G_1(\U_n,r,c)=\emptyset$.
Hence the desired result follows.
$\blacksquare$


There is also a stochastic ordering between $\g_{{}_{n,2}}(F,r,c)$ and $\g_{{}_{n,2}}(\U,r,c)$
provided that $F$ satisfies some regularity conditions.

\begin{proposition}
\label{prop:stoch-order}
Let $\X_n=\{X_1,X_2,\ldots,X_n\}$ be a random sample from
an absolutely continuous distribution $F$ with $\mS(F)\subseteq(0,1)$.
If
\begin{equation}
\label{eqn:stoch-order}
F\bigl( X_{(n)}/r \bigr)< F\left(X_{(n)}\right)/r \text{ and }
F\bigl( X_{(1)} \bigr) < r\,F\left( \left( X_{(1)}+r-1 \right)/r \right)+1-r \text{ hold a.s., }
\end{equation}
then $\g_{{}_{n,2}}(F,r,c)<^{ST}\g_{{}_{n,2}}(\U,r,F(c))$
where $<^{ST}$ stands for ``stochastically smaller than".
If $<$'s in \eqref{eqn:stoch-order} are replaced with $>$'s,
then $\g_{{}_{n,2}}(\U,r,F(c))<^{ST}\g_{{}_{n,2}}(F,r,c)$.
If $<$'s in expression \eqref{eqn:stoch-order} are replaced with $=$'s,
then $\g_{{}_{n,2}}(F,r,c)\stackrel{d}{=}\g_{{}_{n,2}}(\U,r,F(c))$ where $\stackrel{d}{=}$ stands
for equality in distribution.
\end{proposition}

\noindent {\bfseries Proof:}
Let $U_i$ and 
$U_{(i)}$ be as in proof of Proposition \ref{prop:NF vs NPE}. 
Then the parameter $c$ for $N(\cdot,r,c)$ with $\X_n$ in $(0,1)$
corresponds to $F(c)$ for $\U_n$.
Then the $\G_1$-region for $\U_n$ based on $N(\cdot,r,F(c))$ is
$\G_1(\U_n,r,F(c))=(U_{(n)}/r,F(c) ] \cup [F(c),\left(U_{(1)}+r-1\right)/r )$;
likewise, $\G_1(\X_n,r,c)=(X_{(n)}/r,M_c ] \cup [M_c,\left(X_{(1)}+r-1\right)/r )$.
But the conditions in \eqref{eqn:stoch-order} imply that
$\G_1(\U_n,r,F(c)) \subsetneq F(\G_1(\X_n,r,c))$,
since such an $F$ preserves order.
So $\U_n \cap F(\G_1(\X_n,r,c)) = \emptyset$ implies that
$\U_n \cap \G_1(\U_n,r,F(c)) = \emptyset$ and
$\U_n \cap F(\G_1(\X_n,r,c)) = \emptyset$ iff $\X_n \cap \G_1(\X_n,r,F(c)) = \emptyset$.
Hence
$$p_{{}_n}(F,r,c)=P(\X_n \cap \G_1(\X_n,r,c) = \emptyset)<
P(\U_n \cap \G_1(\U_n,r,F(c)) = \emptyset)=p_{{}_n}(\U,r,F(c)).$$
Then $\g_{{}_{n,2}}(F,r,c)<^{ST}\g_{{}_{n,2}}(\U,r,F(c))$ follows.
The other cases follow similarly.
$\blacksquare$


\begin{remark}
We can also find the exact distribution of $\g_{{}_{n,2}}(F,r,c)$
for $F$ whose pdf is piecewise constant
with support in $(0,1)$ as in \cite{ceyhan:dom-num-CCCD-NonUnif}.
Note that the simplest of such distributions is the uniform distribution $\U(0,1)$.
The exact distribution of $\g_{{}_{n,2}}(F,r,c)$ for (piecewise) polynomial
$f(x)$ with at least one piece is of degree 1 or higher and support in $(0,1)$
can be obtained using the multinomial expansion of the term $(\cdot)^{n-2}$
in Equation \eqref{eqn:integrand} with careful bookkeeping.
However, the resulting expression for $p_{{}_n}(F,r,c)$ is extremely lengthy
and not that informative (see \cite{ceyhan:dom-num-CCCD-NonUnif}).

For fixed $n$, one can obtain $p_{{}_n}(F,r,c)$ for $F$
(omitted for the sake of brevity)
by numerical integration of the below expression.
\begin{eqnarray*}
p_{{}_n}(F,r,c)=P\bigl( \g_{{}_{n,2}}(F,r,c)=2 \bigr)&=&
\int\int_{\mS(F)\setminus(\delta_1,\delta_2)}H(x_1,x_n)\,dx_ndx_1,
\end{eqnarray*}
where $H(x_1,x_n)$ is given in Equation $\eqref{eqn:integrand}$. $\square$
\end{remark}

Recall the $\F(\R^d)$-random $\D_{n,m}(r,c)$-digraphs.
We call the digraph which obtains
in the special case of $\Y_2=\{\y_1,\y_2\}$ and support of $F_X$ in $(\y_1,\y_2)$,
\emph{$\F(\y_1,\y_2)$-random $\D_{n,2}(r,c)$-digraph}.
Below, we provide asymptotic results
pertaining to the distribution of such digraphs.

\subsection{The Asymptotic Distribution of the Domination Number of
$\F(\y_1,\y_2)$-random $\D_{n,2}(r,c)$-digraphs}
\label{sec:asy-dist-generalF}
Although the exact distribution of $\g_{{}_{n,2}}(F,r,c)$ is not analytically available
in a simple closed form for $F$
whose density is not piecewise constant,
the asymptotic distribution of $\g_{{}_{n,2}}(F,r,c)$ is available for larger families of distributions.
First, we present the asymptotic distribution of $\g_{{}_{n,2}}(F,r,c)$ for
$\D_{n,2}(r,c)$-digraphs with $\Y_2=\{\y_1,\y_2\} \subset \R$ with $\y_1<\y_2$
for general $F$ with support $\mS(F) \subseteq (\y_1,\y_2)$.
Then we will extend this to the case with $\Y_m \subset \R$ with $m>2$.

Let $c \in (0,1/2)$ and $r \in (1,2)$.
Then for $c=(r-1)/r$, i.e., $M_c=\y_1+(r-1)(\y_2-\y_1)/r$,
we define the family of distributions
\begin{equation*}
\F_1\bigl(\y_1,\y_2\bigr) :=
\Bigl \{\text{$F$ :
$(\y_1,\y_1+\ve) \cup \bigl( M_c,M_c+\ve \bigr) \subseteq \mS(F)\subseteq(\y_1,\y_2)$
for some $\ve \in (0,c)$ with $c=(r-1)/r$} \Bigr\}.
\end{equation*}
Similarly,
let $c \in (1/2,1)$ and $r \in (1,2)$.
Then for $c=1/r$, i.e., $M_c=\y_1+(\y_2-\y_1)/r$ with $r \in (1,2)$,
we define
\begin{equation*}
\F_2\bigl(\y_1,\y_2\bigr) :=
\Bigl \{\text{$F$ :
$(\y_2-\ve,\y_2)\cup \bigl(M_c-\ve,M_c\bigr) \subseteq \mS(F)\subseteq(\y_1,\y_2)$
for some $\ve \in (0,1-c)$ with $c=1/r$} \Bigr\}.
\end{equation*}

Let the $k^{th}$ order right (directed) derivative at $x$ be defined as
$f^{(k)}(x^+):=\lim_{h \rightarrow 0^+}\frac{f^{(k-1)}(x+h)-f^{(k-1)}(x)}{h}$
for all $k \ge 1$ and the right limit at $u$ be defined as $f(u^+):=\lim_{h \rightarrow 0^+}f(u+h)$.
The left derivatives and limits are defined similarly with $+$'s being replaced by $-$'s.

\begin{theorem}
\label{thm:kth-order-gen-(r-1)/r}
\textbf{Main Result 3:}
Let $\Y_2=\{\y_1,\y_2\} \subset \R$ with $-\infty < \y_1 < \y_2<\infty$,
$\X_n=\{X_1,X_2,\ldots,X_n\}$ with $X_i \stackrel {iid}{\sim} F \in \F_1(\y_1,\y_2)$,
and $c \in (0,1/2)$.
Let $D_{n,2}$ be the $\F_1(\y_1,\y_2)$-random $\D_{n,2}(r,c)$-digraph
based on $\X_n$ and $\Y_2$.
\begin{itemize}
\item[(i)]
Then for $n>1$, $r \in (1,\infty)$, and $c=(r-1)/r$
we have $\g_{{}_{n,2}}(F,r,(r-1)/r) \sim 1+ \BER\bigl(p_{{}_n}(F,r,(r-1)/r)\bigr)$.
Note also that $\g_{{}_{1,2}}(F,r,(r-1)/r)=1$ for all $r \ge 1$;
for $r=1$,
we have $\g_{{}_{n,2}}(F,1,0)=1$ for all $n \ge 1$ and
for $r = \infty$,
we have $\g_{{}_{n,2}}(F,\infty,1)=1$ for all $n \ge 1$.

\item[(ii)]
Furthermore,
suppose $k \ge 0$ is the smallest integer for which
$F(\cdot)$ has continuous right derivatives up to order $(k+1)$ at $\y_1$,
$\y_1+(r-1)(\y_2-\y_1)/r$, and
$f^{(k)}(\y_1^+)+r^{-(k+1)}\,f^{(k)}\left( \left( (r-1)(\y_2-\y_1)/r \right)^+ \right) \not= 0$
and $f^{(i)}(\y_1^+)=f^{(i)}\left( \left( \y_1+(r-1)(\y_2-\y_1)/r \right)^+ \right)=0$ for all $i=0,1,2,\ldots,(k-1)$
and suppose also that $F(\cdot)$ has a continuous left derivative at $\y_2$.
Then for bounded $f^{(k)}(\cdot)$, $c=(r-1)/r$, and $r \in (1,2)$,
we have the following limit
$$p(F,r,(r-1)/r) =
\lim_{n \rightarrow \infty}p_{{}_n}(F,r,(r-1)/r) =
\frac{f^{(k)}(\y_1^+)}
{f^{(k)}(\y_1^+)+r^{-(k+1)}\,
f^{(k)}\left( \left( \y_1+(r-1)(\y_2-\y_1)/r \right)^+ \right)}.$$
\end{itemize}
\end{theorem}

Note that in Theorem \ref{thm:kth-order-gen-(r-1)/r}
\begin{itemize}
\item
with $(\y_1,\y_2)=(0,1)$,
we have
$p(F,r,(r-1)/r) =
\frac{f^{(k)}(0^+)}
{f^{(k)}(0^+)+r^{-(k+1)}\,f^{(k)}\left( \left( (r-1)/r \right)^+ \right)}$,
\item
if $f^{(k)}(\y_1^+)=0$ and
$ f^{(k)} \left( \left( \y_1+(r-1)(\y_2-\y_1)/r \right)^+ \right)\not=0$,
then $p_{{}_n}(F,r,(r-1)/r)\rightarrow 0$ as $n \rightarrow \infty$,
at rate $O\bigl( \kappa_1(f)\cdot n^{-(k+2)/(k+1)} \bigr)$
where $\kappa_1(f)$ is a constant depending on $f$
and
\item
if $f^{(k)}(\y_1^+)\not=0$ and
$f^{(k)} \left( \left( \y_1+(r-1)(\y_2-\y_1)/r \right)^+ \right)=0$,
then $p_{{}_n}(F,r,(r-1)/r) \rightarrow 1$
as $n \rightarrow \infty$, at rate $O\bigl(\kappa_1(f)\cdot n^{-(k+2)/(k+1)} \bigr)$.
\end{itemize}

\begin{theorem}
\label{thm:kth-order-gen-1/r}
\textbf{Main Result 4:}
Let $\Y_2=\{\y_1,\y_2\} \subset \R$ with $-\infty < \y_1 < \y_2<\infty$,
$\X_n=\{X_1,X_2,\ldots,X_n\}$ with $X_i \stackrel {iid}{\sim} F \in \F_2(\y_1,\y_2)$,
and $c \in (1/2,1)$.
Let $D_{n,2}$ be the $\F_2(\y_1,\y_2)$-random $\D_{n,2}(r,c)$-digraph
based on $\X_n$ and $\Y_2$.

\begin{itemize}
\item[(i)]
Then for $n>1$, $r \in (1,\infty)$, and $c=1/r$
we have $\g_{{}_{n,2}}(F,r,1/r) \sim 1+ \BER\bigl(p_{{}_n}(F,r,1/r)\bigr)$.
Note also that $\g_{{}_{1,2}}(F,r,1/r)=1$ for all $r \ge 1$;
for $r=1$,
we have $\g_{{}_{n,2}}(F,1,1)=1$ for all $n \ge 1$ and
for $r = \infty$,
we have $\g_{{}_{n,2}}(F,\infty,0)=1$ for all $n \ge 1$.

\item[(ii)]
Furthermore,
suppose $\ell \ge 0$ is the smallest integer for which
$F(\cdot)$ has continuous left derivatives up to order $(\ell+1)$ at $\y_2$, and $\y_1+(\y_2-\y_1)/r$,
and
$f^{(\ell)}(\y_2^-)+r^{-(\ell+1)}\,f^{(\ell)}\left( \left( \y_1+(\y_2-\y_1)/r \right)^- \right) \not= 0$
and $f^{(i)}(\y_2^-)=f^{(i)}\left( \left( \y_1+(\y_2-\y_1)/r \right)^- \right)=0$ for all $i=0,1,2,\ldots,(\ell-1)$
and suppose also that $F(\cdot)$ has a continuous right derivative at $\y_1$.
Additionally,
for bounded $f^{(\ell)}(\cdot)$,
$c=1/r$, and $r \in (1,2)$
we have the following limit
$$
p(F,r,1/r) =
\lim_{n \rightarrow \infty}p_{{}_n}(F,r,1/r) =
\frac{f^{(\ell)}(\y_2^-)}
{f^{(\ell)}(\y_2^-)+r^{-(\ell+1)}\,f^{(\ell)}\left( \left( \y_1+(\y_2-\y_1)/r \right)^- \right)}.$$
\end{itemize}
\end{theorem}

%
%

Note that in Theorem \ref{thm:kth-order-gen-1/r}
\begin{itemize}
\item
with $(\y_1,\y_2)=(0,1)$,
we have
$p(F,r,1/r) =
\frac{f^{(\ell)}(1^-)}
{f^{(\ell)}(1^-)+r^{-(\ell+1)}\,f^{(\ell)}\left( \left( 1/r \right)^- \right)}$,
\item
if $f^{(\ell)}(\y_2^-)=0$ and
$f^{(\ell)} \left( \left( \y_1+(\y_2-\y_1)/r \right)^- \right) \not=0$,
then $p_{{}_n}(F,r,1/r)\rightarrow 0$ as $n \rightarrow \infty$,
at rate $O\bigl( \kappa_2(f)\cdot n^{-(\ell+2)/(\ell+1)} \bigr)$
where $\kappa_2(f)$ is a constant depending on $f$
and
\item
if $f^{(\ell)}(\y_2^-)\not=0$ and
$f^{(\ell)} \left( \left( \y_1+(\y_2-\y_1)/r \right)^- \right)=0$,
then $p_{{}_n}(F,r,1/r) \rightarrow 1$
as $n \rightarrow \infty$, at rate $O\bigl(\kappa_2(f)\cdot n^{-(\ell+2)/(\ell+1)} \bigr)$.
\end{itemize}

\begin{remark}
\label{rem:unbounded}
In Theorems \ref{thm:kth-order-gen-(r-1)/r} and \ref{thm:kth-order-gen-1/r},
we assume that $f^{(k)}(\cdot)$ and $f^{(\ell)}(\cdot)$
are bounded on $(\y_1,\y_2)$, respectively.
If $f^{(k)}(\cdot)$ is not bounded on $(\y_1,\y_2)$ for $k \ge 0$,
in particular at $\y_1$, and $\y_1+(r-1)(\y_2-\y_1)/r$, for example,
$\lim_{x \rightarrow \y_1^+}f^{(k)}(x)=\infty$,
then we have
$$p(F,r,(r-1)/r) =\lim_{\delta \rightarrow 0^+}\frac{f^{(k)}(\y_1+\delta)}
{\left[f^{(k)}(\y_1+\delta)+r^{-(k+1)}\,f^{(k)} \left( (\y_1+(r-1)(\y_2-\y_1)/r)+\delta \right)\right]}.$$
If $f^{(\ell)}(\cdot)$ is not bounded on $(\y_1,\y_2)$ for $\ell \ge 0$,
in particular at $\y_1+(\y_2-\y_1)/r$, and $\y_2$, for example,
$\lim_{x \rightarrow \y_2^-}f^{(\ell)}(x)=\infty$,
then we have
$$p(F,r,1/r) =\lim_{\delta \rightarrow 0^+}\frac{f^{(\ell)}(\y_2-\delta)}
{\left[f^{(\ell)}(\y_2-\delta)+r^{-(\ell+1)}\,f^{(\ell)} \left( (\y_1+(\y_2-\y_1)/r)-\delta \right)\right]}.\;\;\square$$
\end{remark}

\begin{remark}
The rates of convergence in  Theorems \ref{thm:kth-order-gen-(r-1)/r} and
\ref{thm:kth-order-gen-1/r} depends on $f$.
From the proofs of Theorems \ref{thm:kth-order-gen-(r-1)/r} and \ref{thm:kth-order-gen-1/r},
it follows that for sufficiently large $n$,
$$p_n(F,r,(r-1)/r) \approx p(F,r,(r-1)/r) +\frac{\kappa_1(f)}{n^{-(k+2)/(k+1)}}
\text{ and }
p_n(F,r,1/r) \approx p(F,r,1/r) +\frac{\kappa_2(f)}{n^{-(\ell+2)/(\ell+1)}},$$
where
$\kappa_1(f)=\frac{s_1\,s_3^{\frac{1}{k+1}}+s_2\,\G \left(\frac{k+2}{k+1} \right)}
{(k+1)\,s_3^{\frac{k+2}{k+1}} }$
with
$\G(x)=\int_{0}^{\infty} e^{-t}t^{(x-1)}\, dt$,
$s_1=\frac{1}{n^{k+1}k!}\,f^{(k)}(\y_1^+)$,
$s_2=\frac{1}{n(k+1)!}\, f^{(k+1)}(\y_1^+)$,
and
$s_3=\frac{1}{(k+1)!}p(F,r,(r-1)/r)$,
$\kappa_2(f)=\frac{q_1\,\G\left( \frac{\ell+2}{\ell+1} \right)+ q_2\,q_3^{\frac{1}{\ell+1}}}
{(\ell+1)\,q_3^{\frac{\ell+2}{\ell+1}}}$,
$q_1=\frac{(-1)^{\ell+1}}{n(\ell+1)!}\,f^{(\ell+1)}(\y_2^-)$,
$q_2=\frac{(-1)^{\ell}}{n^{\ell+1}\ell!}\,f^{(\ell)}(\y_2^-)$,
and
$q_3=\frac{(-1)^{\ell+1}}{(\ell+1)!}p(F,r,1/r)$
provided the derivatives exist. $\square$
\end{remark}

The conditions of the Theorems \ref{thm:kth-order-gen-(r-1)/r} and \ref{thm:kth-order-gen-1/r}
might seem a bit esoteric. 
However, most of the well known functions that are scaled
and properly transformed to be pdf of some random variable
with support in $(\y_1,\y_2)$ satisfy the conditions
for some $k$ or $\ell$,
hence one can compute the corresponding limiting probabilities
$p(F,r,(r-1)/r)$ and $p(F,r,1/r)$.

\begin{example}
(a)
For example, with $F=\U(\y_1,\y_2)$, in Theorem \ref{thm:kth-order-gen-(r-1)/r},
we have
$k=0$ and $f(\y_1^+)=f\left( (\y_1+(r-1)(\y_2-\y_1)/r)^+ \right)=1/(\y_2-\y_1)$,
and
in Theorem \ref{thm:kth-order-gen-1/r}, we have
$\ell=0$ and $f(\y_2^-)=f\left( (\y_1+(\y_2-\y_1)/r)^- \right)=1/(\y_2-\y_1)$.
Then $\lim _{n \rightarrow \infty}p_n(F,r,(r-1)/r)=\lim _{n \rightarrow \infty}p_n(F,r,1/r)=r/(r+1)$,
which agrees with the result given in Equation \eqref{eqn:asymptotic-uniform}.

(b)
For $F$ with pdf $f(x)=\bigl( x+1/2 \bigr)\,\I\bigl( 0 <x<1 \bigr)$,
we have $k=0$, $f(0^+)=1/2$, and $f\left( (\frac{r-1}{r})^+ \right)=3/2-1/r$ in Theorem \ref{thm:kth-order-gen-(r-1)/r}.
Then $p(F,r,(r-1)/r)=\frac{r^2}{r^2+3r-2}$.
As for Theorem \ref{thm:kth-order-gen-1/r},
we have $\ell=0$,  $f(1^-)=3/2$ and $f\left( (\frac{1}{r})^- \right)=1/r+1/2$.
Then $p(F,r,1/r)=\frac{3r^2}{3r^2+r+2}$.

(c)
For $F$ with pdf
$f(x)=(\pi/2)|\sin(2 \pi x)|\I(0 < x < 1)=
(\pi/2)(\sin(2 \pi x)\I(0 < x \le 1/2)-\sin(2 \pi x)\I(1/2 < x <1 ))$,
we have $k=0$, $f(0^+)=0$,
and $f\left( (\frac{r-1}{r})^+ \right)=(\pi/2)(\sin(2 \pi (r-1)/r)$ in Theorem \ref{thm:kth-order-gen-(r-1)/r}.
Then $p(F,r,(r-1)/r)=0$.
As for Theorem \ref{thm:kth-order-gen-1/r},
we have $\ell=0$,
$f(1^-)=0$ and $f\left( (\frac{1}{r})^- \right)=-(\pi/2)(\sin(2 \pi/r)$.
Then $p(F,r,1/r)=0$.

(d)
For $F$ with pdf
$f(x)=\frac{\pi\,r}{4(r-1)}\sin(\pi\,r x/(r-1))\I(0 < x \le (r-1)/r)+ g(x)\I((r-1)/r < x <1 )$,
where $g(x)$ is a nonnegative function such that
$\int_{(r-1)/r}^1 g(t)dt=1/2$,
we have $k=1$, $f'(0^+)=\frac{(\pi r)^2}{4(r-1)^2}$,
and $f'\left( (\frac{r-1}{r})^+ \right)=\frac{(\pi r)^2}{4(r-1)^2}$ in Theorem \ref{thm:kth-order-gen-(r-1)/r}.
Then $p(F,r,(r-1)/r)=r^2/(r^2-1)$.

(e)
For the beta distribution with parameters $\nu_{1,n},\nu_{2,n}$,
denoted $\beta(\nu_{1,n},\nu_{2,n})$,
where $\nu_{1,n}, \nu_{2,n} \ge 1$,
the pdf is given by
$$f(x,{\nu_{1,n}},{\nu_{2,n}})= {\frac{x^{\nu_{1,n}-1}(1-x)^{\nu_{2,n}-1}}
{\beta(\nu_{1,n},\nu_{2,n}) }} \;\I(0<x<1) \text{ where }\beta(\nu_{1,n},\nu_{2,n})=\frac{\G(\nu_{1,n})\,\G(\nu_{2,n})}{ \G(\nu_{1,n} +\nu_{2,n})}.$$
Then
in Theorem \ref{thm:kth-order-gen-(r-1)/r}
we have $k=0$, $f(0^+)=0$,
and $f\left( (\frac{r-1}{r})^+ \right)=\frac{(r-1)^{\nu_{1,n}-1}}{r^{\nu_{1,n}+\nu_{2,n}-1}}$.
So $p(\beta(\nu_{1,n},\nu_{2,n}),r,(r-1)/r)=0$.
As for Theorem \ref{thm:kth-order-gen-1/r},
we have $\ell=0$, $f(1^-)=0$,
and $f\left( (\frac{1}{r})^- \right)=\frac{(r-1)^{\nu_{2,n}-1}}{r^{\nu_{1,n}+\nu_{2,n}-1}}$.
Then $p(\beta(\nu_{1,n},\nu_{2,n}),r,1/r)=0$.

(f)
Consider $F$ with pdf $f(x)=\left(\pi \sqrt{x\,(1-x)}\right)^{-1} \;\I(0<x<1)$.
Notice that $f(x)$ is unbounded at $x \in \{0,1\}$.
Using Remark \ref{rem:unbounded},
it follows that $p(F,r,(r-1)/r)=p(F,r,1/r)=1$.
$\square$
\end{example}


\begin{remark}
\label{rem:same-asydist-uniform}
In Theorem \ref{thm:kth-order-gen-(r-1)/r},
if we have $f^{(k)}(0^+)=f^{(k)} \left( (\frac{r-1}{r})^+ \right)$,
then
$\lim_{n \rightarrow \infty}p_{{}_n}(F,r,(r-1)/r)=\frac{1}{1+r^{-(k+1)}}.$
In particular, if $k=0$, then
$\lim_{n \rightarrow \infty}p_{{}_n}(F,r,(r-1)/r)=r/(r+1)$.
Hence
$\g_{{}_{n,2}}(F,r,(r-1)/r)$ and $\g_{{}_{n,2}}(\U,r,(r-1)/r)$
have the same asymptotic distribution.

In Theorem \ref{thm:kth-order-gen-1/r},
if we have $f^{(\ell)}(1^-)=f^{(\ell)} \left( (\frac{1}{r})^- \right)$,
then
$\lim_{n \rightarrow \infty}p_{{}_n}(F,r,1/r)=\frac{1}{1+r^{-(\ell+1)}}.$
In particular, if $\ell=0$, then
$\lim_{n \rightarrow \infty}p_{{}_n}(F,r,1/r)=r/(r+1)$.
Hence
$\g_{{}_{n,2}}(F,r,1/r)$ and $\g_{{}_{n,2}}(\U,r,1/r)$
have the same asymptotic distribution.
$\square$
\end{remark}

The asymptotic distribution of  $\g_{{}_{n,2}}(F,r,c)$ for $r=2$ and $c=1/2$
is as follows (see \cite{ceyhan:dom-num-CCCD-NonUnif} for its derivation).
\begin{theorem}
\label{thm:kth-order-gen}
Let
$\F\bigl(\y_1,\y_2\bigr):=\Bigl \{F:
(\y_1,\y_1+\ve) \cup (\y_2-\ve,\y_2)\cup \bigl( (\y_1+\y_2)/2-\ve,(\y_1+\y_2)/2+\ve \bigr)
\subseteq \mS(F)
\subseteq(\y_1,\y_2) \text{ for some } \ve \in (0,(\y_1+\y_2)/2) \Bigr\}.
$
Let $\Y_2=\{\y_1,\y_2\} \subset \R$ with $-\infty < \y_1 < \y_2<\infty$,
$\X_n=\{X_1,\ldots,X_n\}$ with $X_i \stackrel {iid}{\sim} F \in \F(\y_1,\y_2)$,
and $D_{n,2}$ be the random $\D_{n,2}$-digraph based on $\X_n$ and $\Y_2$.
\begin{itemize}
\item[(i)]
Then for $n>1$, we have
$\g_{{}_{n,2}}(F,2,1/2) \sim 1+ \BER\bigl(p_{{}_n}(F,2,1/2)\bigr)$.
Note also that $\g_{{}_{1,2}}(F,2,1/2)=1$.

\item[(ii)]
Furthermore,
suppose $k \ge 0$ is the smallest integer for which
$F(\cdot)$ has continuous right derivatives up to order $(k+1)$ at $\y_1,\,(\y_1+\y_2)/2$,
$f^{(k)}(\y_1^+)+2^{-(k+1)}\,f^{(k)}\left( \left( \frac{\y_1+\y_2}{2} \right)^+ \right) \not= 0$
and $f^{(i)}(\y_1^+)=f^{(i)} \left( \left( \frac{\y_1+\y_2}{2} \right)^+ \right)=0$ for all $i=0,1,\ldots,k-1$;
and $\ell \ge 0$ is the smallest integer for which
$F(\cdot)$ has continuous left derivatives up to order $(\ell+1)$ at $\y_2,\,(\y_1+\y_2)/2$,
$f^{(\ell)}(\y_2^-)+2^{-(\ell+1)}\,f^{(\ell)}\left( \left( \frac{\y_1+\y_2}{2} \right)^- \right) \not= 0$
and $f^{(i)}(\y_2^-)=f^{(i)}\left( \left( \frac{\y_1+\y_2}{2} \right)^- \right)=0$ for all $i=0,1,\ldots,\ell-1$.
Additionally,
for bounded $f^{(k)}(\cdot)$ and $f^{(\ell)}(\cdot)$,
we have the following limit
$$
p(F,2,1/2) =
\lim_{n \rightarrow \infty}p_{{}_n}(F,2,1/2) =
\frac{f^{(k)}(\y_1^+)\,f^{(\ell)}(\y_2^-)}
{\left[f^{(k)}(\y_1^+)+2^{-(k+1)}\,f^{(k)} \left( \left( \frac{\y_1+\y_2}{2} \right)^+ \right)\right]\,\left[f^{(\ell)}(\y_2^-)+
2^{-(\ell+1)}\,f^{(\ell)}\left( \left( \frac{\y_1+\y_2}{2} \right)^- \right)\right]}.$$
\end{itemize}
\end{theorem}

Notice the interesting behavior of $p(F,r,c)$ around $(r,c)=(2,1/2)$.
There is a jump (hence discontinuity) in $p(F,r,(r-1)/r)$ and in $p(F,r,1/r)$ at $r=2$.

\section{The Distribution of the Domination Number of $\D_{n,m}(r,c)$-digraphs}
\label{sec:dist-multiple-intervals}
We now consider the more challenging case of $m>2$.
For $\omega_1<\omega_2$ in $\R$, define the family of distributions
$$
\mathscr H(\R):=\bigl \{ F_{X,Y}:\;(X_i,Y_i) \sim F_{X,Y} \text{ with support }
\mS(F_{X,Y})=(\omega_1,\omega_2)^2 \subsetneq \R^2,\;\;X_i \sim F_X \text{ and } Y_i \stackrel{iid}{\sim}F_Y \bigr\}.
$$
We provide the exact distribution of $\g_{{}_{n,m}}(F,r,c)$ for
$\mathscr H(\R)$-random digraphs in the following theorem.
Let $[m]:=\bigl\{ 0,1,2,\ldots,m-1 \bigr\}$ and
$\Theta^S_{a,b}:=\bigl\{ (u_1,u_2,\ldots u_b):\;\sum_{i=1}^{b}u_i = a:\; u_i \in S, \;\;\forall i \bigr\}$.
If $Y_i$ have a continuous distribution,
then the order statistics of $\Y_m$ are distinct a.s.
Given $Y_{(i)}=\y_{(i)}$ for $i=1,2,\ldots,m$,
let
$\vec{n}$ be the vector of numbers $n_i$,
$f_{\vec{Y}}(\vec{\y})$ be the joint distribution of the order statistics of $\Y_m$,
i.e., $f_{\vec{Y}}(\vec{\y})=\frac{1}{m!}\prod_{i=1}^m f(\y_i)\,\I(\omega_1<\y_1<\y_2<\ldots<\y_m<\omega_2)$,
and $f_{i,j}(\y_i,\y_j)$ be the joint distribution of $Y_{(i)},Y_{(j)}$.
Then we have the following theorem.

\begin{theorem}
\label{thm:general-Dnm}
Let $D$ be an $\mathscr H(\R)$-random $\D_{n,m}(r,c)$-digraph
with $n>1$, $m>1$, $r \in [1,\infty)$ and $c \in (0,1)$.
Then the probability mass function of the domination number of D is given by
{\small
$$P(\g_{{}_{n,m}}(F,r,(r-1)/r)=q)=\int_{\mathscr S} \sum_{\vec{n} \in \Theta^{[n+1]}_{n,(m+1)}}
\sum_{\vec{q}\in \Theta^{[3]}_{q,(m+1)}} P(\vec{N}=\vec{n})\,\zeta(q_1,n_1)\,\zeta(q_{m+1},\,n_{m+1})
\prod_{j=2}^{m}\eta(q_i,n_i)f_{\vec{Y}}(\vec{\y})\,d\y_1 \ldots d\y_m$$
}
where
$P(\vec{N}=\vec{n})$ is the joint probability of $n_i$ points
falling into intervals $\mI_i$ for $i=0,1,2,\ldots,m$,
$q_i \in \{0,1,2\}$, $q=\sum_{i=0}^m q_i$ and
\begin{align*}
\zeta(q_i,n_i)&=\max\bigl( \I(n_i=q_i=0),\I(n_i \ge q_i=1) \bigr) \text{ for } i=1,(m+1),
\text{ and }\\
\eta(q_i,n_i)&=\max \bigl( \I(n_i=q_i=0),\I(n_i \ge q_i \ge 1) \bigr)\cdot
p_{{}_{n_i}}(F_i,r,(r-1)/r))^{\I(q_i=2)}\,\bigl( 1-p_{{}_{n_i}}(F_i,r,(r-1)/r)) \bigr)^{\I(q_i=1)}\\
& \text{ for $i=1,2,3,\ldots,(m-1),$ and the region of integration is given by}\\
\mathscr S:=\bigl\{&(\y_1,\y_2,\ldots,\y_m)\in
(\omega_1,\omega_2)^2:\,\omega_1<\y_1<\y_2<\ldots<\y_m<\omega_2 \bigr\}.
\end{align*}
The special cases of $n=1$, $m=1$, $r \in \{1,\infty\}$ and $c \in \{0,1\}$ are
as in Theorem \ref{thm:gamma-Dnm-r-M}.
\end{theorem}
Proof is as in Theorem 6.1 of \cite{ceyhan:dom-num-CCCD-NonUnif}.
A similar construction is available for $c=1/r$.

This exact distribution for finite $n$ and $m$ has a simpler form
when $\X$ and $\Y$ points are both uniform in a bounded interval in $\R$.
Define $\mathscr U(\R)$ as follows
$$
\mathscr U(\R):=\bigl \{ F_{X,Y}: \text{ $X$ and $Y$ are independent}
\;X_i \stackrel{iid}{\sim} \U(\omega_1,\omega_2) \text{ and } Y_i \stackrel{iid}{\sim}\U(\omega_1,\omega_2),
\text{ with } -\infty <\omega_1<\omega_2<\infty \bigr\}.
$$
Clearly, $\mathscr U(\R) \subsetneq \mathscr H(\R)$.
Then we have the following corollary to Theorem \ref{thm:general-Dnm}.
\begin{corollary}
\label{cor:uniform-Dnm}
Let $D$ be a $\mathscr U(\R)$-random $\D_{n,m}(r,c)$-digraph
with $n>1$, $m>1$, $r \in [1,\infty)$ and $c \in (0,1)$.
Then the probability mass function of the domination number of $D$ is given by
$$P(\g_{{}_{n,m}}(r,(r-1)/r)=q)=\frac{n!m!}{(n+m)!}
\sum_{\vec{n} \in \Theta^{[n+1]}_{n,(m+1)}}
\sum_{\vec{q}\in \Theta^{[3]}_{q,(m+1)}}
\zeta(q_1,n_1)\,\zeta(q_{m+1},\,n_{m+1})
\prod_{j=2}^{m}\eta(q_i,n_i).$$
The special cases of $n=1$, $m=1$, $r \in \{1,\infty\}$ and $c \in \{0,1\}$ are
as in Theorem \ref{thm:gamma-Dnm-r-M}.
\end{corollary}
Proof is as in Theorem 2 of \cite{priebe:2001}.
A similar construction is available for $c=1/r$.
For $n,m < \infty$, the expected value of domination number is
\begin{equation}
\label{eqn:expected-gamma-Dnm}
\E[\g_{{}_{n,m}}(F,r,c)]=P\left(X_{(1)}<Y_{(1)}\right)+P\left(X_{(n)} > Y_{(m)}\right)+
\sum_{i=1}^{m-1}\sum_{k=1}^n\,P(N_i=k)\,\E[\g_{{}_{n_i,2}}(F_i,r,c)]
\end{equation}
where
\begin{multline*}
P(N_i=k)=\\
\int_{\omega_1}^{\omega_2}\int_{\y_{(i)}}^{\omega_2}
f_{i-1,i}\left(\y_{(i)},\y_{(i+1)}\right) \Bigl[F_X\left(\y_{(i+1)}\right)-
F_X\left(\y_{(i)}\right)\Bigr]^k\Bigl[1-\left(F_X\left(\y_{(i+1)}\right)-
F_X\left(\y_{(i)}\right)\right)\Bigr]^{n-k}\,d\y_{(i+1)}d\y_{(i)}
\end{multline*}
and $\E[\g_{{}_{n_i,2}}(F_i,r,c)]=1+p_n(F_i,r,c)$.
Then as in Corollary 6.2 of \cite{ceyhan:dom-num-CCCD-NonUnif},
we have
\begin{corollary}
\label{cor:Egnn goes infty}
For $F_{X,Y} \in \mathscr H(\R)$ with support $\mS(F_X) \cap \mS(F_Y)$ of positive measure
with $r \in [1,\infty)$ and $c \in (0,1)$,
we have $\lim_{n \rightarrow \infty}\E[\g_{{}_{n,n}}(F,r,c)] = \infty$.
\end{corollary}

\begin{theorem}
\label{thm:asy-general-Dnm}
\textbf{Main Result 5:}
Let $D_{n,m}(r,c)$ be an $\mathscr H(\R)$-random $\D_{n,m}(r,c)$-digraph.
Then
\begin{itemize}
\item[(i)]
for fixed $n<\infty$, $\lim_{m \rightarrow \infty}\g_{{}_{n,m}}(F,r,c)=n$ a.s. for all $r \ge 1$ and $c \in [0,1]$.
\item[] For fixed $m<\infty$, and
\item[(ii)]
for $r=1$ and $c \in (0,1)$,
$\lim_{n \rightarrow \infty}P(\g_{{}_{n,m}}(F,1,c)=2m)=1$
and
$\lim_{n \rightarrow \infty}P(\g_{{}_{n,m}}(F,1,0)=m+1)=\\
\lim_{n \rightarrow \infty}P(\g_{{}_{n,m}}(F,1,1)=m+1)=1$
\item[(iii)]
for $r>2$ and $c \in (0,1)$, $\lim_{n \rightarrow \infty}P(\g_{{}_{n,m}}(F,r,c)=m+1)=1$,
\item[(iv)]
for $r=2$,
if $c \not= 1/2$,
then $\lim_{n \rightarrow \infty}P(\g_{{}_{n,m}}(F,2,c)=m+1)=1$;\\
if $c = 1/2$,
then $\lim_{n \rightarrow \infty}\g_{{}_{n,m}}(F,2,1/2)\stackrel{d}{=}m+1+\BIN(m,p(F_i,2,1/2))$,
\item[(v)]
for $r \in (1,2)$,
if $r \not= \tau=\max(c,1-c)$,
then $\lim_{n \rightarrow \infty} \g_{{}_{n,m}}(F,r,c)$ is degenerate;
otherwise,
it is non-degenerate.
That is, for $r \in [1,2)$,
as $n \rightarrow \infty$,
\begin{equation}
\label{eqn:asy-unif-rM-mult-int}
\g_{{}_{n,m}}(F,r,c) \sim
\left\lbrace \begin{array}{ll}
       m+1+\BIN(m,p(F_i,r,c)),  & \text{for $r = 1/\tau$,}\\
       m+1,           & \text{for $r > 1/\tau$,}\\
       2m+1,           & \text{for $r < 1/\tau$.}\\
\end{array} \right.
\end{equation}
\end{itemize}
\end{theorem}

\begin{itemize}
\item[] {\bfseries Proof:}
\item[]
Part (i) is trivial.
Part (ii) follows from Theorems \ref{thm:gamma-Dnm-r-M} and \ref{thm:gamma-Dnm-r=1-M},
since as $n_i \rightarrow \infty$,
we have
$\X_{[i]} \not= \emptyset$ a.s. for all $i$.
\item[]
Part (iii) follows from Theorem \ref{thm:r and M-asy},
since for $c \in (0,1)$,
it follows that $r > 1/\tau $ implies $r>2$
and as $n_i \rightarrow \infty$,
we have
$\g_{{}_{n_i,2}}(F_i,r,c) \rightarrow 1$ in probability for all $i$.
\item[]
In part (iv), for $r=2$ and $c \not= 1/2$,
based on Theorem \ref{thm:r=2 and M_c=c},
as $n_i \rightarrow \infty$,
we have
$\g_{{}_{n_i,2}}(F_i,r,c) \rightarrow 1$ in probability for all $i$.
The result for $r=2$ and $c = 1/2$ is proved in \cite{ceyhan:dom-num-CCCD-NonUnif}.
\item[]
Part (v) follows from Theorem \ref{thm:r and M-asy}. $\blacksquare$
\end{itemize}

\begin{remark}
\label{rem:ext-multi-dim}
\textbf{Extension to Higher Dimensions:}

Let $\Y_m=\left \{\y_1,\y_2,\ldots,\y_m \right\}$ be $m$ points in
general position in $\R^d$ and $T_i$ be the $i^{th}$ Delaunay cell
in the Delaunay tessellation (assumed to exist) based on $\Y_m$
for $i=1,2,\ldots,J_m$. 
Let $\X_n$ be a random sample from a distribution $F$ in
$\R^d$ with support $\mS(F) \subseteq \C_H(\Y_m)$
where $\C_H(\Y_m)$ stands for the convex hull of $\Y_m$.
In $\R$
a Delaunay tessellation is an intervalization (i.e.,
partitioning of $\R$ by intervals),
provided that no two points in $\Y_m$ are concurrent. 

We define the proportional-edge proximity region in $\R^2$.
The extension to $\R^d$ with $d>2$ is straightforward
(see \cite{ceyhan:dom-num-NPE-MASA} for an explicit extension).
Let $\TY$ be the triangle (including the interior) with vertices $\Y_3=\{\y_1,\y_2,\y_3\}$,
$e_i$ be the edge opposite vertex $\y_i$,
and $M_i$ be the midpoint of edge $e_i$ for $i=1,2,3$.
We first construct the vertex regions based on a point $M \in \R^2 \setminus \Y_3$
called \emph{$M$-vertex regions},
by the lines joining $M$ to a point on each of the edges of $\TY$.
See \cite{ceyhan:Phd-thesis} for a more general definition of vertex regions.
Preferably, $M$ is selected to be in the interior of the triangle $\TY$.
For such an $M$, the corresponding vertex regions can be defined using
a line segment joining $M$ to $e_i$. 
With center of mass $M_{CM}$,
the lines joining $M_{CM}$ and $\Y_3$ are the \emph{median lines}
which cross edges at midpoints $M_i$ for $i=1,2,3$.
The vertex regions in Figure \ref{fig:prox-map-def} are
center of mass vertex regions. 
For $r \in [1,\infty]$, define $N(\cdot,r,M)$
to be the \emph{(parametrized) proportional-edge proximity map} with
$M$-vertex regions as follows
(see also Figure \ref{fig:prox-map-def} with $M=M_{CM}$ and $r=2$).
Let $R_M(v)$ be the vertex region associated with vertex $v$ and $M$.
For $x \in \TY \setminus \Y_3$,
let $v(x) \in \Y_3$ be the vertex whose region contains $x$; i.e., $x \in R_M(v(x))$.
If $x$ falls on the boundary of two $M$-vertex regions,
we assign $v(x)$ arbitrarily.
Let $e(x)$ be the edge of $\TY$ opposite of $v(x)$,
$\ell(v(x),x)$ be the line parallel to $e(x)$ and passes through $x$,
and $d(v(x),\ell(v(x),x))$ be the Euclidean
distance from $v(x)$ to $\ell(v(x),x)$.
For $r \in [1,\infty)$,
let $\ell_r(v(x),x)$ be the line parallel to $e(x)$ such that
$$
d(v(x),\ell_r(v(x),x)) = r\,d(v(x),\ell(v(x),x))\\
\text{ and }\\
d(\ell(v(x),x),\ell_r(v(x),x)) < d(v(x),\ell_r(v(x),x)).
$$
Let $T_r(x)$ be the triangle similar to and with the same
orientation as $\TY$ having $v(x)$ as a vertex and $\ell_r(v(x),x)$
as the opposite edge.
Then the \emph{proportional-edge proximity region} $\N(x,r,M)$ is defined to be $T_r(x) \cap \TY$.
Notice that $\ell(v(x),x)$ divides the two edges of $T_r(x)$
(other than the one lies on $\ell_r(v(x),x)$) proportionally with the factor $r$.
Hence the name \emph{proportional-edge proximity region}.

\begin{figure} [ht]
\centering
\scalebox{.35}{\input{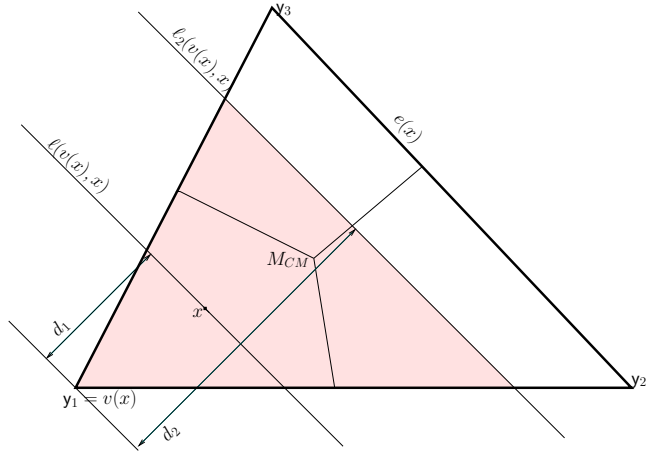}}
\caption{
\label{fig:prox-map-def}
Construction of proportional-edge proximity region, $\N(x,2,M_{CM})$ (shaded region)
for an $x$ in the CM-vertex region for $\y_1$, $R_{CM}(\y_1)$
where $d_1=d(v(x),\ell(v(x),x))$ and
$d_2=d(v(x),\ell_2(v(x),x))=2\,d(v(x),\ell(v(x),x))$.
}
\end{figure}

Notice that in $\R$, $M$ is the center parametrized by $c$,
e.g., the center of mass $M_{CM}$ corresponds to $c=1/2$,
but for other $M \in \TY$, there is no direct counterpart in $\R$.
The vertex regions in $\R$ with $\Y_2=\{\y_1,\y_2\}$ are $(\y_1,M_c)$ and $(M_c,\y_2)$.
Observe that $N(x,r,c)$ in $\R$ is an open interval,
while in $\R^d$, the region $\N(x,r,M)$ is a closed region.
However, the interiors of $\N(x,r,M)$ satisfy the class cover problem of \cite{cannon:2000}.
$\square$
\end{remark}

\section{Discussion}
\label{sec:disc-conclusions}
In this article,
we present the distribution of the domination number
of a random digraph family called proportional-edge proximity catch digraph (PCD)
which is based on two classes of points.
Points from one of the classes 
constitute the vertices of the PCDs,
while the points from the other class
are used in the binary relation
that assigns the arcs of the PCDs.

We introduce the proximity map which is the one-dimensional version
of $N(\cdot,r,c)$ of \cite{ceyhan:dom-num-NPE-SPL} and \cite{ceyhan:dom-num-NPE-MASA}.
This proximity map can also be viewed as an extension of the proximity map
of \cite{priebe:2001} and \cite{ceyhan:dom-num-CCCD-NonUnif}.
The PCD we consider is based on a parametrized proximity map
in which there is an expansion parameter $r$ and a centrality parameter $c$.
We provide the exact and asymptotic distributions of
the domination number for proportional-edge PCDs
for uniform data
and compute the asymptotic distribution for non-uniform data
for the entire range of $(r,c)$.
The results in this article can also be
viewed as generalizations of the main results of
\cite{priebe:2001} and \cite{ceyhan:dom-num-CCCD-NonUnif} in several directions.
\cite{priebe:2001} provided the exact distribution of the
domination number of class cover catch digraphs (CCCDs) based on $\X_n$ and $\Y_m$
both of which were sets of iid random variables from
uniform distribution on $(\omega_1,\omega_2) \subset \R$ with $-\infty<\omega_1<\omega_2<\infty$
and the proximity map $N(x):=B(x,r(x))$ where $r(x):=\min_{\y \in \Y_m}d(x,\y)$.
\cite{ceyhan:dom-num-CCCD-NonUnif} investigates the
distribution of the domination number of CCCDs for non-uniform data
and provides the asymptotic distribution for a large family
of (non-uniform) distributions.
%
Furthermore, this article will form the foundation of the generalizations and calculations
for uniform and non-uniform cases in multiple dimensions.
As in \cite{ceyhan:dom-num-NPE-SPL},
we can use the domination number in testing
one-dimensional spatial point patterns
and our results will help make the power comparisons
possible for a large family of distributions.

We demonstrate an interesting behavior of the domination number of
proportional-edge PCD for one-dimensional data.
For uniform data or data from a distribution which satisfies some regularity conditions
(see Section \ref{sec:asy-dist-generalF})
and fixed $1<n<\infty$,
the distribution of the domination number is a translated form of
(extended) binomial distribution
$\BIN(m,p_{{}_n}(F,r,c))$ where $m$ is the number of (inner) intervals
and $p_{{}_n}(F,r,c)$ is the probability that the domination number
of the proportional-edge PCD is 2.
Here $p_{{}_n}(F,r,c)$ is allowed to take 0 or 1 also.
For finite $n>1$,
the parameter, $p_{{}_n}(\U,r,c)$, of distribution of the domination number
of the proportional-edge PCD based on uniform data is continuous in $r$ and $c$
for all $r \ge 1$ and $c \in (0,1)$
and has jumps (hence discontinuities) at $c=0,1$.
For fixed $(r,c) \in [1,\infty) \times (0,1)$,
the parameter, $p(\U,r,c)$, of the asymptotic distribution
exhibits some discontinuities.
For $c \in (0,1)$
the asymptotic distribution is nondegenerate at $r=1/\max(c,1-c)$.
The asymptotic distribution of the domination number is degenerate for the expansion parameter $r > 2$.
For $r \in (1,2]$, there exists threshold values for the centrality parameter $c$ for which
the asymptotic distribution is non-degenerate.
In particular, for $c \in \{ (r-1)/r,1/r \}$ with $r \in (1,2]$
the asymptotic distribution
of the domination number is a translated form of $\BIN(m,p(\U,r,c)$
where
$p(F,r,c)$ is continuous in $r$.
Additionally,
by symmetry,
we have $p(\U,r,(r-1)/r)=p(\U,r,1/r)$ for $r \in (1,2]$.
For $r > 1/\max(c,1-c)$
the domination number converges in probability to 1,
and for $r < 1/\max(c,1-c)$
the domination number converges in probability to 2.
On the other hand,
at $(r,c)= (2,1/2)$,
the asymptotic distribution is again a translated form of $\BIN(m,p(\U,2,1/2))$
but there is yet another jump at $r=2$,
as $p(\U,2,1/2)=4/9$ while $\lim_{r \rightarrow 2} p(\U,r,(r-1)/r)=\lim_{r \rightarrow 2} p(\U,r,1/r)=2/3$.
This second jump might be due to the symmetry
for the domination number at $c=1/2$. 

\section*{Acknowledgments}
Supported by TUBITAK Kariyer Project Grant 107T647.


\section*{APPENDIX}
In the proofs of Theorems \ref{thm:r=2 and M_c=c},
\ref{thm:r and M_c=1/2},
\ref{thm:r and M},
based on Proposition \ref{prop:scale-inv-NYr},
we can assume $(\y_1,\y_2)=(0,1)$.

\subsection*{Proof of Theorem \ref{thm:r=2 and M_c=c} }
Given $X_{(1)}=x_1$ and $X_{(n)}=x_n$,
let $a=x_n/2$ and $b=(1+x_1)/2$
and due to symmetry, we only consider $c  \in (0,1/2)$.
For $r=2$,
the $\G_1$-region is $\G_1(\X_n,2,c)=(a,c] \cup [c,b)$
so we have two cases for $\G_1$-region:
case (1) $\G_1(\X_n,2,c)=(a,b)$ which occurs when $a<c<b$,
and
case (2) $\G_1(\X_n,2,c)=[c,b)$ which occurs when $c<a<b$.
The cases in which $b<a$ and $b<c$ are not possible,
since $c<1/2$, $b>1/2$, and $a<1/2$.

\noindent \textbf{Case (1)} $\G_1(\X_n,2,c)=(a,b)$, i.e., $a<c<b$:
For $c \in [1/3,1/2]$,
we have
\begin{multline}
\label{eqn:r=2-c-in-[1/3,1/2)-case1}
P(\g_{{}_{n,2}}(\U,2,c)=2,\; \G_1(\X_n,2,c)=(a,b),\; c \in [1/3,1/2) )=\\
\left( \int_{0}^{1/3}\int_{(1+x_1)/2}^{2\,c}+\int_{1/3}^{c}\int_{2\,x_1}^{2\,c} \right)
n(n-1)f(x_1)f(x_n)\bigl[F(x_n)-F(x_1)+F(a)-F(b)\bigr]^{(n-2)}\,dx_ndx_1=\\
\frac{4}{9}\, \left( 3\,c-\frac{1}{2} \right)^n-{\frac{8}{9}}\,{4}^{-n}-{\frac{8}{9}}\, \left( \frac{3\,c-1}{2} \right)^n.
\end{multline}
For $c \in (0,1/3)$,
we have
\begin{multline}
\label{eqn:r=2-c-in-(0,1/3)-case1}
P(\g_{{}_{n,2}}(\U,2,c)=2,\; \G_1(\X_n,2,c)=(a,b),\; c \in (0,1/3))=\\
\int_{0}^{4c-1}\int_{(1+x_1)/2}^{2\,c}
n(n-1)f(x_1)f(x_n)\bigl[F(x_n)-F(x_1)+F(a)-F(b)\bigr]^{(n-2)}\,dx_ndx_1=\\
\frac{4}{9}\, \left( 1 -3\,c \right)^n+\frac{4}{9}\, \left( 3\,c-\frac{1}{2} \right)^n-{
\frac{8}{9}}\,{4}^{-n}.
\end{multline}

\noindent \textbf{Case (2)} $\G_1(\X_n,2,c)=[c,b)$, i.e., $c<a<b$:
For $c \in [1/3,1/2)$,
we have
\begin{multline}
\label{eqn:r=2-c-in-[1/3,1/2)-case2}
P(\g_{{}_{n,2}}(\U,2,c)=2,\;\G_1(\X_n,2,c)=[c,b), \;c \in [1/3,1/2))=\\
\int_{0}^{4c-1}\int_{(1+x_1)/2}^{2\,c}
n(n-1)f(x_1)f(x_n)\bigl[F(x_n)-F(x_1)+F(c)-F(b)\bigr]^{(n-2)}\,dx_ndx_1=\\
\frac{2}{3}\, \left( \frac{3\,c-1}{2} \right)^n-\frac{2}{3}\, \left( \frac{1-c}{2} \right)^n
-\frac{2}{3}\, \left( 3\,c-\frac{1}{2} \right)^n+\frac{2}{3}\, \left( c + \frac{1}{2} \right)^n.
\end{multline}

For $c \in (0,1/3)$,
we have
\begin{multline}
\label{eqn:r=2-c-in-(0,1/3)-case2}
P(\g_{{}_{n,2}}(\U,2,c)=2,\;\G_1(\X_n,2,c)=[c,b), \;c \in (0,1/3))=\\
\left( \int_{0}^{4c-1}\int_{2c}^{1}+\int_{4c-1}^{c}\int_{(1+x_1)/2}^{1} \right)
n(n-1)f(x_1)f(x_n)\bigl[F(x_n)-F(x_1)+F(c)-F(b)\bigr]^{(n-2)}\,dx_ndx_1=\\
\frac{2}{3} \left(c+\frac{1}{2}\right)^n-\frac{1}{3} (1-3\,c)^n-\frac{2}{3}\left(3c -\frac{1}{2}\right)^n-\frac{2}{3}\left(\frac{1-c}{2}\right)^n
\end{multline}
The probability
$P(\g_{{}_{n,2}}(\U,2,c)=2,\;c \in [1/3,1/2))$ is the sum of probabilities
in \eqref{eqn:r=2-c-in-[1/3,1/2)-case1} and \eqref{eqn:r=2-c-in-[1/3,1/2)-case2},
and
$P(\g_{{}_{n,2}}(\U,2,c)=2,\;c \in (0,1/3))$ is the sum of probabilities
in \eqref{eqn:r=2-c-in-(0,1/3)-case1} and \eqref{eqn:r=2-c-in-(0,1/3)-case2}.
By symmetry,
$P(\g_{{}_{n,2}}(\U,2,c)=2,\;c \in [1/2,1))=P(\g_{{}_{n,2}}(\U,2,1-c)=2,\; 1-c \in (0,1/2]$.
The special cases for $c=\{0,1\}$ follow by construction.
$\blacksquare$

\subsection*{Proof of Theorem \ref{thm:r and M_c=1/2} }

Given $X_{(1)}=x_1$ and $X_{(n)}=x_n$, let $a=x_n/r$ and $b=(x_1+r-1)/r$.
For $r \ge 1$ and $c=1/2$,
the $\G_1$-region is $\G_1(\X_n,r,1/2)=(a,1/2] \cup [1/2,b)$,
so we have four cases for $\G_1$-region:
case (1) $\G_1(\X_n,r,1/2)=(a,b)$ which occurs when $a<1/2<b$,
case (2) $\G_1(\X_n,r,1/2)=(a,1/2]$ which occurs when $a<b<1/2$ or $b<a<1/2$,
case (3) $\G_1(\X_n,r,1/2)=[1/2,b)$ which occurs when $1/2<a<b$ or $1/2<b<a$,
and
case (4) $\G_1(\X_n,r,1/2)=\emptyset$ which occurs when $b<1/2<a$.
Cases (2) and (3) are symmetric, so they yield the same probabilities.\\
\noindent \textbf{Case (1)} $\G_1(\X_n,r,1/2)=(a,b)$, i.e., $a<1/2<b$:
For $r \ge 2$,
we have
\begin{multline}
\label{eqn:c=1/2-r-in-[2,infty)-case1}
P(\g_{{}_{n,2}}(\U,r,1/2)=2,\;\G_1(\X_n,r,1/2)=(a,b),\;r \ge 2)=\\
\left( \int_{0}^{1/(r+1)}\int_{(x_1+r-1)/r}^{1}+\int_{1/(r+1)}^{1/r}\int_{r\,x_1}^{1} \right)
n(n-1)f(x_1)f(x_n)\bigl[F(x_n)-F(x_1)+F(a)-F(b)\bigr]^{(n-2)}\,dx_ndx_1=\\
\frac{2\,r}{(r+1)^2} \left[\left(\frac{2}{r}\right)^{n-1}-\left(\frac{r-1}{r^2}\right)^{n-1}\frac{1}{r}\right].
\end{multline}
For $1 \le r < 2$,
we have
\begin{multline}
\label{eqn:c=1/2-r-in-[1,2)-case1}
P(\g_{{}_{n,2}}(\U,r,1/2)=2,\;\G_1(\X_n,r,1/2)=(a,b),\;1 \le r < 2)=\\
\left( \int_{1-r/2}^{1/(r+1)}\int_{(x_1+r-1)/r}^{r/2}+\int_{1/(r+1)}^{1/2}\int_{r\,x_1}^{r/2} \right)
n(n-1)f(x_1)f(x_n)\bigl[F(x_n)-F(x_1)+F(a)-F(b)\bigr]^{(n-2)}\,dx_ndx_1=\\
\frac{r^2(r-1)^n}{(r+1)^2} \left[1-\left(\frac{r-1}{2\,r}\right)^{n-1}\right]
.
\end{multline}

\noindent \textbf{Case (2)} $\G_1(\X_n,r,1/2)=(a,1/2]$, i.e., $a<b<1/2$ or $b<a<1/2$:
Here, $r \ge 2$ is not possible,
since $x_1 < 1-r/2$.
For $1 \le r < 2$,
we have
\begin{multline}
\label{eqn:c=1/2-r-in-[1,2)-case2}
P(\g_{{}_{n,2}}(\U,r,1/2)=2,\;\G_1(\X_n,r,1/2)=(a,1/2],\;1 \le r < 2)=\\
\int_{0}^{1-r/2}\int_{1/2}^{r/2}
n(n-1)f(x_1)f(x_n)\bigl[F(x_n)-F(x_1)+F(a)-F(1/2)\bigr]^{(n-2)}\,dx_ndx_1=\\
\frac{r}{(r+1)}\left[\left(\frac{r}{2}\right)^n-(r-1)^n-
\left(\frac{1}{2\,r}\right)^{n}+\left(\frac{(r-1)^2}{2\,r}\right)^n\right].
\end{multline}
By symmetry, \noindent \textbf{Case (3)} yields the same result as \noindent \textbf{Case (2)}.

\noindent \textbf{Case (4)} $\G_1(\X_n,r,1/2)=\emptyset$, i.e., $b<1/2<a$:
Here, $r \ge 2$ is not possible, since $x_n-x_1>r-1$.
For $1 \le r < 2$,
we have
\begin{multline}
\label{eqn:c=1/2-r-in-[1,2)-case4}
P(\g_{{}_{n,2}}(\U,r,1/2)=2,\;\G_1(\X_n,r,1/2)=\emptyset,\;1 < r < 2)=\\
\int_{0}^{1-r/2}\int_{1/2}^{1}
n(n-1)f(x_1)f(x_n)\bigl[F(x_n)-F(x_1)\bigr]^{(n-2)}\,dx_ndx_1=
1+(r-1)^n-2\left(\frac{r}{2}\right)^n.
\end{multline}

The probability
$P(\g_{{}_{n,2}}(\U,r,1/2)=2,\;r \ge 2)$ is as in \eqref{eqn:c=1/2-r-in-[2,infty)-case1},
and
$P(\g_{{}_{n,2}}(\U,r,1/2)=2,\;r \in [1,\infty))$
is the sum of probabilities in \eqref{eqn:c=1/2-r-in-[1,2)-case1},
\eqref{eqn:c=1/2-r-in-[1,2)-case4},
and twice the probability in \eqref{eqn:c=1/2-r-in-[1,2)-case2}. $\blacksquare$

\subsection*{Proof of Theorem \ref{thm:r and M} }

Given $X_{(1)}=x_1$ and $X_{(n)}=x_n$,
let $a=x_n/r$ and $b=(x_1+r-1)/r$
and assume $c \in (0,1/2)$.
There are two cases for $c$, namely,
\textbf{Case I-} $c \in ((3-\sqrt{5})/2,1/2)$
and
\textbf{Case II-} $c \in (0,(3-\sqrt{5})/2]$

\textbf{Case I-}
For $r \ge 1$ and $c \in ((3-\sqrt{5})/2,1/2)$,
the $\G_1$-region is $\G_1(\X_n,r,c)=(a,c] \cup [c,b)$.
So
we have four cases for $\G_1$-region:
case (1) $\G_1(\X_n,r,c)=(a,b)$ which occurs when $a<c<b$,
case (2) $\G_1(\X_n,r,c)=(a,c]$ which occurs when $a<b<c$ or $b<a<c$,
case (3) $\G_1(\X_n,r,c)=[c,b)$ which occurs when $c<a<b$ or $c<b<a$, and
case (4) $\G_1(\X_n,r,c)=\emptyset$ which occurs when $b<c<a$.\\
\noindent \textbf{Case (1)}
$\G_1(\X_n,r,c)=(a,b)$, i.e., $a<c<b$:
In this case,
for $r \ge 1/c$,
we have
\begin{multline}
\label{eqn:c-in-I-r-in-[1/c,infty)-case1}
P(\g_{{}_{n,2}}(\U,r,c)=2,\;\G_1(\X_n,r,c)=(a,b),\;r \ge 1/c)=\\
\left( \int_{0}^{1/(r+1)}\int_{(x_1+r-1)/r}^{1}+\int_{1/(r+1)}^{1/r}\int_{r\,x_1}^{1} \right)
n(n-1)f(x_1)f(x_n)\bigl[F(x_n)-F(x_1)+F(a)-F(b)\bigr]^{(n-2)}\,dx_ndx_1=\\
\frac{2\,r}{(r+1)^2} \left( \left(\frac{2}{r}\right)^{n-1}-
\left(\frac{r-1}{r^2} \right)^{n-1} \right)
\end{multline}
For $1/(1-c) \le r < 1/c$,
we have
\begin{multline}
\label{eqn:c-in-I-r-in-[1/(1-c),1/c)-case1}
P(\g_{{}_{n,2}}(\U,r,c)=2,\;\G_1(\X_n,r,c)=(a,b),\;1/(1-c) \le r < 1/c)=\\
\left( \int_{0}^{1/(r+1)}\int_{(x_1+r-1)/r}^{c\,r}+\int_{1/(r+1)}^{c}\int_{r\,x_1}^{c\,r} \right)
n(n-1)f(x_1)f(x_n)\bigl[F(x_n)-F(x_1)+F(a)-F(b)\bigr]^{(n-2)}\,dx_ndx_1=\\
\frac{r^2}{(r+1)^2} \Biggl[ \left(c(r+1)-\frac{r-1}{r}\right)^n
- \left(\frac{r-1}{r}\right)^{n-1} \left((c\,r+c-1)^n-\frac{1}{r^n}\right) \Biggr].
\end{multline}
For $(1-c)/c \le r < 1/(1-c)$,
we have
\begin{multline}
\label{eqn:c-in-I-r-in-[(1-c)/c,1/(1-c))-case1}
P(\g_{{}_{n,2}}(\U,r,c)=2,\;\G_1(\X_n,r,c)=(a,b),\;(1-c)/c \le r < 1/(1-c))=\\
\left( \int_{r(c-1)+1}^{1/(r+1)}\int_{(x_1+r-1)/r}^{c\,r}+\int_{1/(r+1)}^{c}\int_{r\,x_1}^{c\,r} \right)
n(n-1)f(x_1)f(x_n)\bigl[F(x_n)-F(x_1)+F(a)-F(b)\bigr]^{(n-2)}\,dx_ndx_1=\\
\frac{r^2(r-1)^{n-1}}{(r+1)^2}\left[ (r-1)-\frac{1}{r^{n-1}}[(r-c\,r-c)^n-(c\,r+c-1)^n] \right].
\end{multline}

For $1 \le r < (1-c)/c$,
we have
\begin{multline}
\label{eqn:c-in-I-r-in-[1,(1-c)/c)-case1}
P(\g_{{}_{n,2}}(\U,r,c)=2,\;\G_1(\X_n,r,c)=(a,b),\;1 \le r < (1-c)/c)=\\
\int_{r(c-1)+1}^{c\,r^2-r+1}\int_{(x_1+r-1)/r}^{c\,r}
n(n-1)f(x_1)f(x_n)\bigl[F(x_n)-F(x_1)+F(a)-F(b)\bigr]^{(n-2)}\,dx_ndx_1=\\
\frac{r^2 (r-1)^{n-1}}{(r+1)^2}
\left[r-1+(1-c\,r-c)^n+ \frac{(r-c\,r-c)^n}{r^{n-1}} \right].
\end{multline}

\noindent \textbf{Case (2)}
$\G_1(\X_n,r,c)=(a,c]$, i.e., $a<b<c$ or $b<a<c$:
Here $r \ge 1/(1-c)$ is not possible, since $r(c-1)+1>0$.
Hence $r \ge 1/c$ is not possible either.
For $1 \le r < 1/(1-c)$,
we have
\begin{multline}
\label{eqn:c-in-I-r-in-[1,1/(1-c))-case2}
P(\g_{{}_{n,2}}(\U,r,c)=2,\;\G_1(\X_n,r,c)=(a,c],\;1 \le r < 1/(1-c))=\\
\int_{0}^{r(c-1)+1}\int_{c}^{c\,r}
n(n-1)f(x_1)f(x_n)[F(x_n)-F(x_1)+F(a)-F(c)]^{(n-2)}\,dx_ndx_1=\\
\frac{r}{r+1} \left[c^n\left(r^n-\frac{1}{r^n}\right)-(r-1)^n \left(1-\frac{r-c\,r-c}{r}\right)^n \right].
\end{multline}

\noindent \textbf{Case (3)}
$\G_1(\X_n,r,c)=[c,b)$, i.e., $c<a<b$ or $c<b<a$:
Here $r \ge 1/c$ is not possible, since $x_n >c\,r$.
For $1/(1-c) \le r < 1/c$,
we have
\begin{multline}
\label{eqn:c-in-I-r-in-[1/(1-c),1/c)-case3}
P(\g_{{}_{n,2}}(\U,r,c)=2,\;\G_1(\X_n,r,c)=[c,b),\;1/(1-c) \le r < 1/c)=\\
\int_{0}^{c}\int_{c\,r}^{1}
n(n-1)f(x_1)f(x_n)\bigl[F(x_n)-F(x_1)+F(c)-F(b)\bigr]^{(n-2)}\,dx_ndx_1=\\
\frac{1}{(r+1)r^{n-1}} \left[
(r-1)^n (c\,r-1+c)^n
+ (1+c\,r)^n - (c\,r^2+c\,r-r+1)^n - (1-c)^n
\right].
\end{multline}
For $1 \le r < 1/(1-c)$,
we have
\begin{multline}
\label{eqn:c-in-I-r-in-[1,1/(1-c))-case3}
P(\g_{{}_{n,2}}(\U,r,c)=2,\;\G_1(\X_n,r,c)=[c,b),\;1 \le r < 1/(1-c))=\\
\int_{r(c-1)+1}^{c}\int_{c\,r}^{1}
n(n-1)f(x_1)f(x_n)\bigl[F(x_n)-F(x_1)+F(c)-F(b)\bigr]^{(n-2)}\,dx_ndx_1=\\
\frac{r}{r+1} \left[
(r-1)^n \left(\left(\frac{c\,r-1+c}{r}\right)^n-1\right)+
(1-c)^n\left(r^n-\frac{1}{r^n}\right)
\right].
\end{multline}

\noindent \textbf{Case (4)}
$\G_1(\X_n,r,c)=\emptyset$, i.e., $b<c<a$:
Here, $r \ge 1/c$ is not possible, since $x_n>c\,r$;
and $1/(1-c) \le r < 1/c$ is not possible, since $x_1 < r(c-1)+1$.
For $1 \le r < 1/(1-c)$,
we have
\begin{multline}
\label{eqn:c-in-I-r-in-[1,1/(1-c))-case4}
P(\g_{{}_{n,2}}(\U,r,c)=2,\;\G_1(\X_n,r,c)=\emptyset,\;1 \le r < 1/(1-c))=\\
\int_{0}^{r(c-1)+1}\int_{c\,r}^{1}
n(n-1)f(x_1)f(x_n)(F(x_n)-F(x_1))^{(n-2)}\,dx_ndx_1=
1+ (r-1)^n - r^n [c^n + (1-c)^n].
\end{multline}

The probability $P(\g_{{}_{n,2}}(\U,r,c)=2,\;r \ge 1/c)$
is the same as in \eqref{eqn:c-in-I-r-in-[1/c,infty)-case1};
$P(\g_{{}_{n,2}}(\U,r,c)=2,\;1/(1-c) \le r < 1/c)$
is the sum of probabilities in \eqref{eqn:c-in-I-r-in-[1/(1-c),1/c)-case1}
and \eqref{eqn:c-in-I-r-in-[1/(1-c),1/c)-case3};
$P(\g_{{}_{n,2}}(\U,r,c)=2,\;(1-c)/c \le r < 1/(1-c))$
is the sum of probabilities in \eqref{eqn:c-in-I-r-in-[(1-c)/c,1/(1-c))-case1},
\eqref{eqn:c-in-I-r-in-[1,1/(1-c))-case2},
\eqref{eqn:c-in-I-r-in-[1,1/(1-c))-case3},
and
\eqref{eqn:c-in-I-r-in-[1,1/(1-c))-case4};
$P(\g_{{}_{n,2}}(\U,r,c)=2,\;1 \le r < (1-c)/c)$
is the sum of probabilities in \eqref{eqn:c-in-I-r-in-[1,(1-c)/c)-case1},
\eqref{eqn:c-in-I-r-in-[1,1/(1-c))-case2},
\eqref{eqn:c-in-I-r-in-[1,1/(1-c))-case3},
and
\eqref{eqn:c-in-I-r-in-[1,1/(1-c))-case4}.

\textbf{Case II-}
For $r \ge 1$ and $c \in (0,(3-\sqrt{5})/2)$,
we have the same cases for the $\G_1$-region is as above.\\
\textbf{Case (1)}:
For $r \ge 1/c$,
the probability
$P(\g_{{}_{n,2}}(\U,r,c)=2,\;\G_1(\X_n,r,c)=(a,b),\;r \ge 1/c)$ as in
\eqref{eqn:c-in-I-r-in-[1/c,infty)-case1}.

For $(1-c)/c \le r < 1/c$,
the probability
$P(\g_{{}_{n,2}}(\U,r,c)=2,\;\G_1(\X_n,r,c)=(a,b),\;(1-c)/c \le r < 1/c)$ is as in
\eqref{eqn:c-in-I-r-in-[1/(1-c),1/c)-case1}.

For $1/(1-c) \le r < (1-c)/c$,
we have
\begin{multline}
\label{eqn:c-in-II-r-in-[(1/1-c),(1-c)/c)-case1}
P(\g_{{}_{n,2}}(\U,r,c)=2,\;\G_1(\X_n,r,c)=(a,b),\;1/(1-c) \le r < (1-c)/c)=\\
\int_{0}^{c\,r^2-r+1}\int_{(x_1+r-1)/r}^{c\,r}
n(n-1)f(x_1)f(x_n)\bigl[F(x_n)-F(x_1)+F(a)-F(b)\bigr]^{(n-2)}\,dx_ndx_1=\\
\frac{r^2}{(r+1)^2}
\left[
(r-1)^{n-1}\left((1-c\,r-c)^n-\frac{1}{r^{2n-1}}\right)+
\left(\frac{c\,r^2+c\,r-r+1}{r}\right)^n
\right].
\end{multline}

For $1 \le r < 1/(1-c)$,
the probability
$P(\g_{{}_{n,2}}(\U,r,c)=2,\;\G_1(\X_n,r,c)=(a,b),\;1 \le r < 1/(1-c))$ is as in
\eqref{eqn:c-in-I-r-in-[1,(1-c)/c)-case1}.

\textbf{Cases (2)} and \textbf{(4)} are as before.

\noindent \textbf{Case (3)}:
For $(1-c)/c \le r < 1/c$,
the probability
$P(\g_{{}_{n,2}}(\U,r,c)=2,\;\G_1(\X_n,r,c)=[c,b),\;(1-c)/c \le r < 1/c)$ is as in
\eqref{eqn:c-in-I-r-in-[1/(1-c),1/c)-case3}.

For $1/(1-c) \le r < (1-c)/c$,
the probability
$P(\g_{{}_{n,2}}(\U,r,c)=2,\;\G_1(\X_n,r,c)=[c,b),\;1/(1-c) \le r < (1-c)/c)$ is as in
\eqref{eqn:c-in-I-r-in-[1/(1-c),1/c)-case3}.

For $1 \le r < 1/(1-c)$,
the probability
$P(\g_{{}_{n,2}}(\U,r,c)=2,\;\G_1(\X_n,r,c)=[c,b),\;1 \le r < 1/(1-c))$ is as in
\eqref{eqn:c-in-I-r-in-[1,1/(1-c))-case3}.

The probability $P(\g_{{}_{n,2}}(\U,r,c)=2,\;r \ge 1/c)$
is the same as in \eqref{eqn:c-in-I-r-in-[1/c,infty)-case1};
$P(\g_{{}_{n,2}}(\U,r,c)=2,\;(1-c)/c \le r < 1/c)$
is the sum of probabilities in \eqref{eqn:c-in-I-r-in-[1/(1-c),1/c)-case1}
and \eqref{eqn:c-in-I-r-in-[1/(1-c),1/c)-case3};
$P(\g_{{}_{n,2}}(\U,r,c)=2,\;1/(1-c) \le r < (1-c)/c)$
is the sum of probabilities in
\eqref{eqn:c-in-I-r-in-[1/(1-c),1/c)-case3}
and
\eqref{eqn:c-in-II-r-in-[(1/1-c),(1-c)/c)-case1};
$P(\g_{{}_{n,2}}(\U,r,c)=2,\;1 \le r < 1/(1-c))$
is the sum of probabilities in \eqref{eqn:c-in-I-r-in-[1,(1-c)/c)-case1},
\eqref{eqn:c-in-I-r-in-[1,1/(1-c))-case2},
\eqref{eqn:c-in-I-r-in-[1,1/(1-c))-case3},
and
\eqref{eqn:c-in-I-r-in-[1,1/(1-c))-case4}.

By symmetry,
$P(\g_{{}_{n,2}}(\U,r,c)=2,\;c \in [1/2,1))=P(\g_{{}_{n,2}}(\U,r,1-c)=2,\; 1-c \in (0,1/2]$.
The special case for $c\in \{0,1\}$ follows trivially by construction.
$\blacksquare$

\subsection*{Proof of Theorem \ref{thm:r and M-asy} }
Let $c \in (0,1/2)$.
Then $\tau=1-c$.
We first consider \textbf{Case I}: $c \in ((3-\sqrt{5})/2,1/2)$.
In Theorem \ref{thm:r and M},
for $r \ge 1/c > 2$,
it follows that
$\lim_{n \rightarrow \infty}p_{{}_n}(\U,r,c)=\lim_{n \rightarrow \infty}\pi_{1,n}(r,c)=0$,
since $2/r <1 $ and $\frac{r-1}{r^2}<1$.

For $1/(1-c) < r < 1/c$,
we have
$\frac{1+c\,r}{r}<1$ (since $r>1/(1-c)$),
$\frac{1-c}{r}<1$,
$\frac{(c\,r^2-r+c\,r+1)}{r}<1$, 
$\frac{r-1}{r^2}<1$ (since $r-1<r<r^2$),
and
$\frac{(r-1)(c\,r-1+c)}{r}<1$.
Hence
for $1/(1-c) < r < 1/c$
$\lim_{n \rightarrow \infty}p_{{}_n}(\U,r,c)=\lim_{n \rightarrow \infty}\pi_{2,n}(r,c)=0$

For $(1-c)/c < r < 1/(1-c)$,
we have
$r-1<1$ (since $r<1/(1-c)<2$),
$\frac{(r-1)(c\,r-1+c)}{r}<1$,
$\frac{(r-1)(r-c\,r-c)}{r}<1$,
$c/r < 1$ (since $c<r$),
$(1-c)/r < 1$ (since $1-c<r$),
$c\,r<1$ and $(1-c)r<1$ (since $r<1/(1-c)<1/c$).
Hence
$\lim_{n \rightarrow \infty}p_{{}_n}(\U,r,c)=\lim_{n \rightarrow \infty}\pi_{3,n}(r,c)=1$.

For $1 \le r < (1-c)/c$,
we have
$r-1<1$ (since $r<1/(1-c)<2$),
$(r-1)(1-c\,r-c)<1$,
$\frac{(r-1)(1-c\,r-c)}{r}<1$,
$\frac{(r-1)(r-c\,r-c)}{r}<1$,
$c/r < 1$ (since $c<r$),
$(1-c)/r < 1$,
$c\,r<1$ and $(1-c)r<1$.
Hence
$\lim_{n \rightarrow \infty}p_{{}_n}(\U,r,c)=\lim_{n \rightarrow \infty}\pi_{4,n}(r,c)=1$,

But for $r=1/(1-c)$ or $c=(r-1)/r$,
we have
\begin{equation}
p_{{}_n}(\U,r,(r-1)/r)=
\frac{r}{(r+1)^2}
\left[
(r+1)-(r-1)^{n-1}\left(\frac{r^2-r-1}{r^2}\right)^n+
(r-1)^n-
\frac{r+1}{r^{2n}}-
\left(\frac{r-1}{r^2}\right)^{n-1}
\right].
\end{equation}
Letting $n \rightarrow \infty$,
we get $p_{{}_n}(\U,r,(r-1)/r) \rightarrow r/(r+1)$ for $r \in (1,2)$,
since
$\frac{r-1}{r^2}<1$,
$r-1<1$,
$\frac{1}{r^2}<1$,
and
$\frac{(r-1)(r^2-r-1)}{r^2}<1$.

Next we consider \textbf{Case II:} $c \in (0,(3-\sqrt{5})/2]$.
In Theorem \ref{thm:r and M},
for $r \ge 1/c > 2$,
it follows that
$\lim_{n \rightarrow \infty}p_{{}_n}(\U,r,c)=\lim_{n \rightarrow \infty}\vartheta_{1,n}(r,c)=0$,
since $2/r <1 $ and $\frac{r-1}{r^2}<1$.

For $(1-c)/c < r < 1/c$,
we have
$\frac{1+c\,r}{r}<1$ (since $r>1/(1-c)$),
$\frac{1-c}{r}<1$,
$\frac{(c\,r^2-r+c\,r+1)}{r}<1$, 
$\frac{r-1}{r}<1$ (since $r-1<r<r^2$),
and
$\frac{(r-1)(c\,r-1+c)}{r}<1$.
Hence
$\lim_{n \rightarrow \infty}p_{{}_n}(\U,r,c)=\lim_{n \rightarrow \infty}\vartheta_{2,n}(r,c)=0$

For $1/(1-c) < r < (1-c)/c$,
we have
$(r-1)(1-c\,r-c)<1$,
$\frac{(r-1)(1-c\,r-c)}{r}<1$,
$\frac{r-1}{r}<1$,
$\frac{r-1}{r^2}<1$,
$(1+cr)/r < 1$,
$\frac{c\,r^2-c+c\,r+1}{r}<1$,
$c/r < 1$ (since $c<r$),
$(1-c)/r < 1$ (since $1-c<r$),
$c\,r<1$ and $(1-c)r<1$ (since $r<1/(1-c)<1/c$).
Hence
$\lim_{n \rightarrow \infty}p_{{}_n}(\U,r,c)=\lim_{n \rightarrow \infty}\vartheta_{3,n}(r,c)=0$.

For $1 \le r < 1/(1-c)$,
we have
$r-1<1$ (since $r<1/(1-c)<2$),
$(r-1)(1-c\,r-c)<1$,
$\frac{(r-1)(1-c\,r-c)}{r}<1$,
$\frac{(r-1)(r-c\,r-c)}{r}<1$,
$c/r < 1$ (since $c<r$),
$(1-c)/r < 1$,
$c\,r<1$ and $(1-c)r<1$.
Hence
$\lim_{n \rightarrow \infty}p_{{}_n}(\U,r,c)=\lim_{n \rightarrow \infty}\vartheta_{4,n}(r,c)=1$,

But for $r=1/(1-c)$ or $c=(r-1)/r$,
we have
\begin{multline}
p_{{}_n}(\U,r,(r-1)/r)=
\frac{r}{(r+1)^2}
\Biggl[
(r+1)+(r+1)\left(\frac{(r-1)(r^2-r-1)}{r^2}\right)^n-
(r-1)^n+\\
(-1)^n\left(\frac{r-1}{r}\right)^{n-1}(r^2-r-1)^n-
\frac{r-1}{r^{2n}}+
\left(\frac{r-1}{r^2}\right)^{n-1}
\Biggr].
\end{multline}
Letting $n \rightarrow \infty$,
we get $p_{{}_n}(\U,r,(r-1)/r) \rightarrow r/(r+1)$ for $r \in (1,2)$,
since
$r-1<1$,
$\frac{1}{r^2}<1$,
$\frac{(r-1)(r^2-r-1)}{r^2}<1$,
and
$r^2-r-1 <1$.

For $c \in (1/2,1)$, we have $\tau=c$.
By symmetry, the above results follow with $c$ being replaced by $1-c$
and
as $n \rightarrow \infty$,
we get $p_{{}_n}(\U,r,1/r) \rightarrow r/(r+1)$ for $r \in (1,2)$.
Hence the desired result follows.
$\blacksquare$

\subsection*{Proof of Theorem \ref{thm:kth-order-gen-(r-1)/r} }

Case (i) follows trivially from Theorem \ref{thm:gamma 1 or 2}.
The special cases for $n=1$ and $r=\{1,\infty\}$ follow by construction.

\noindent
Case (ii):
Suppose $(\y_1,\y_2)=(0,1)$ and $c \in (0,1/2)$.
Recall that
$\G_1(\X_n,r,c)=(X_{(n)}/r, M_c ] \bigcup [M_c, ( X_{(1)}+r-1 )/r ) \subset (0,1)$
and $\g_{{}_{n,2}}(F,r,c)=2 \text{ iff } \X_n \cap \G_1(\X_n,r,c)=\emptyset$.
Then for finite $n$,
\begin{equation*}
p_{{}_n}(F,r,c)=P\bigl( \g_{{}_{n,2}}(F,r,c)=2 \bigr)=\int_{\mS(F) \setminus (\delta_1,\delta_2)} H(x_1,x_n)\,dx_ndx_1,
\end{equation*}
where $(\delta_1,\delta_2)=\G_1(\X_n,r,c)$ and $H(x_1,x_n)$ is as in Equation \eqref{eqn:integrand}.

Let $\ve \in (0,(r-1)/r)$ and $c=(r-1)/r$.
Then $P\bigl( X_{(1)} <\ve,\; X_{(n)} > 1-\ve \bigr)\rightarrow 1$
as $n \rightarrow \infty$ with the rate of convergence depending on $F$.
Moreover, for sufficiently large $n$,
$(X_{(1)}+r-1)/r > (r-1)/r$ a.s.;
in fact, $(X_{(1)}+r-1)/r \downarrow (r-1)/r$ as $n \rightarrow \infty$ (in probability)
and $X_{(n)}/r > \max((r-1)/r,(X_{(1)}+r-1)/r)$ a.s. since $r \in (1,2)$.
Then for sufficiently large $n$,
we have $\G_1(\X_n,r,(r-1)/r)=[(r-1)/r,(X_{(1)}+r-1)/r)$ a.s.
and
\begin{multline}
\label{eqn:asy-g=2-x1xn}
p_{{}_n}(F,r,(r-1)/r) \approx \int_0^{\ve}\int_{1-\ve}^1 n\,(n-1)f(x_1)f(x_n)
\Bigl[F(x_n)-F(x_1)+F\left((r-1)/r \right)-F\left( (x_1+r-1)/r \right) \Bigr]^{n-2}\,dx_ndx_1\\
=\int_0^{\ve} nf(x_1)
\Biggl(\Bigl[1-F(x_1)+F\left((r-1)/r \right)-F\left( (x_1+r-1)/r \right) \Bigr]^{n-1}-\\
\Bigl[1-\varepsilon-F(x_1)+F\left((r-1)/r \right)-F\left( (x_1+r-1)/r \right) \Bigr]^{n-1}\Biggr)\,dx_1\\
\approx
\int_0^{\ve} nf(x_1)
\Bigl[1-F(x_1)+F\left((r-1)/r \right)-F\left( (x_1+r-1)/r \right) \Bigr]^{n-1}\,dx_1.
\end{multline}
 Let
$$G(x_1)=1-F(x_1)+F\left((r-1)/r \right)-F\left( (x_1+r-1)/r \right).$$
The integral in Equation \eqref{eqn:asy-g=2-x1xn} is critical at $x_1=0$,
since $G(0)=1$,
and for $x_1 \in (0,1)$ the
integral converges to 0 as $n \rightarrow \infty$.
Let $\al_i := -\frac{d^{i+1} G(x_1)}{d x_1^{i+1}}\Big|_{(0^+,0^+)}=
f^{(i)}(0^+)+r^{-(i+1)}\,f^{(i)}\left(\left( \frac{r-1}{r} \right)^+\right)$.
Then by the hypothesis of the theorem, we have $\al_i = 0$ and
$f^{(i)}\left(\left( \frac{r-1}{r} \right)^+\right)=0$ for all $i=0,1,2,\ldots,(k-1)$.
So the Taylor series expansions of $f(x_1)$ around $x_1=0^+$ up to order $k$
and $G(x_1)$ around $0^+$ up to order $(k+1)$
so that $x_1 \in (0,\ve)$, are as follows:
$$f(x_1)=\frac{1}{k!}f^{(k)}(0^+)\,x_1^k+O\left( x_1^{k+1} \right)$$
and
$$G(x_1)=G(0^+)+ \frac{1}{(k+1)!}\left(\frac{d^{k+1}G(0^+)}{d x_1^{k+1}}\right)\,x_1^{k+1}
+O\left( x_1^{k+2} \right) \\
=1-\frac{\al_k}{(k+1)!}\,x_1^{k+1}+O\left(x_1^{k+2}\right).$$
Then substituting these expansions in Equation \eqref{eqn:asy-g=2-x1xn},
we obtain
$$
p_{{}_n}(F,r,(r-1)/r) \approx \int_0^{\ve}
n\Biggl[\frac{1}{k!}f^{(k)}(0^+)\,x_1^k+O\left(x_1^{k+1}\right)\Biggr]
\Biggl[1-\frac{\al_k}{(k+1)!}\,x_1^{k+1}+O\left(x_1^{k+2}\right)\Biggr]^{n-1}\,d x_1.
$$

Now we let $x_1=w\,n^{-1/(k+1)}$ to obtain
\begin{multline}
\label{eqn:Pg2-in-(0,1)-(r-1)/r}
p_{{}_n}(F,r,(r-1)/r) \approx
\int_0^{\ve\,n^{1/(k+1)}}n\,
\Biggl[\frac{1}{n^{k/(k+1)}\,k!}f^{(k)}(0^+)w^k+O\left(n^{-1}\right)\Biggr]\\
\Biggl[1-\frac{1}{n}\left(\frac{\al_k}{(k+1)!}\,w^{k+1}+
O\left(n^{-(k+2)/(k+1)}\right)\right)\Biggr]^{n-1}\,\left(\frac{1}{n^{1/(k+1)}}\right)\,dw\\
\text{letting $n \rightarrow
\infty,$~~~~~~~~~~~~~~~~~~~~~~~~~~~~~~~~~~~~~~~~~~~~~~~~~~~~~~~~
~~~~~~~~~~~~~~~~~~~~~~~~~~~~~~~~~~~~~~~~~~~~~~~~~~~~}\\
\approx
\int_0^{\infty} \frac{1}{k!} f^{(k)}(0^+)w^k \,
\exp\left[-\frac{\al_k}{(k+1)!}\,w^{k+1}\right]\, dw
= \frac{f^{(k)}(0^+)}{\al_k}
= \frac{f^{(k)}(0^+)}
{f^{(k)}(0^+)+r^{-(k+1)}\,f^{(k)} \left( \left( \frac{r-1}{r} \right)^+ \right)},
\end{multline}
as $n \rightarrow \infty$ at rate $O(\kappa_1(f)\cdot n^{-(k+2)/(k+1)})$.

For the general case of $\Y=\{\y_1,\y_2\}$, the transformation
$\phi(x)=(x-\y_1)/(\y_2-\y_1)$ maps $(\y_1,\y_2)$ to $(0,1)$ and
the transformed random variables $U=\phi(X_i)$ are distributed with
density $g(u)=(\y_2-\y_1)\,f(\y_1+u(\y_2-\y_1))$ on $(\y_1,\y_2)$.
Replacing $f(x)$ by $g(x)$ in Equation \eqref{eqn:Pg2-in-(0,1)-(r-1)/r},
the desired result follows.
$\blacksquare$

\subsection*{Proof of Theorem \ref{thm:kth-order-gen-1/r} }

Case (i) follows trivially from Theorem \ref{thm:gamma 1 or 2}.
The special cases for $n=1$ and $r=\{1,\infty\}$ follow by construction.

\noindent
Case (ii):
Suppose $(\y_1,\y_2)=(0,1)$ and $c \in (1/2,1)$.
Let $\ve \in (0,1/r)$.
Then $P\bigl( X_{(1)} <\ve,\; X_{(n)} > 1-\ve \bigr)\rightarrow 1$
as $n \rightarrow \infty$ with the rate of convergence depending on $F$.
Moreover, for sufficiently large $n$,
$X_{(n)}/r < 1/r$ a.s.;
in fact, $X_{(n)}/r \uparrow 1/r$ as $n \rightarrow \infty$ (in probability)
and $(X_{(1)}+r-1)/r < \min(1/r,X_{(n)}/r)$ a.s.
Then for sufficiently large $n$, $\G_1(\X_n,r,c)=(X_{(n)}/r,1/r]$
a.s.
and
\begin{multline}
\label{eqn:asy-g=2-x1xn-1/r}
p_{{}_n}(F,r,1/r) \approx \int_{1-\ve}^1 \int_0^{\ve} n\,(n-1)f(x_1)f(x_n)
\Bigl[F(x_n)-F(x_1)+F\left(x_n/r \right)-F\left( 1/r \right) \Bigr]^{n-2}\,dx_1dx_n.\\
=-\int_{1-\ve}^1 n f(x_n)
\Biggl(\Bigl[F(x_n)-F(\varepsilon)+F\left(x_n/r \right)-F\left( 1/r \right) \Bigr]^{n-1}-
\Bigl[F(x_n)+F\left(x_n/r \right)-F\left( 1/r \right) \Bigr]^{n-1}\Biggr)\,dx_n\\
\approx
\int_{1-\ve}^1 n f(x_n)
\Bigl[F(x_n)+F\left(x_n/r \right)-F\left( 1/r \right) \Bigr]^{n-1}\,dx_n.
\end{multline}
 Let
$$G(x_n)=F(x_n)+F\left(x_n/r \right)-F\left( 1/r \right).$$
The integral in Equation \eqref{eqn:asy-g=2-x1xn-1/r} is critical at $x_n=1$,
since $G(1)=1$,
and for $x_n \in (0,1)$ the
integral converges to 0 as $n \rightarrow \infty$.
So we make the change of variables $z_n=1-x_n$,
then $G(x_n)$ becomes
$$G(z_n)=F(1-z_n)+F((1-z_n)/r)-F(1/r),$$
and Equation \eqref{eqn:asy-g=2-x1xn-1/r} becomes
\begin{equation}
\label{eqn:asy-g=2-zn}
p_{{}_n}(F,r,1/r) \approx \int_0^{\ve} n\,f(1-z_n)\left[G(z_n)\right]^{n-1}\,dz_n.
\end{equation}
The new integral is critical at $z_n=0$.
Let
$\be_i := (-1)^{i+1}\frac{d^{i+1} G(z_n)}{d z_n^{i+1}}\Big|_{0^+}=
f^{(i)}(1^-)+r^{-(i+1)}\,f^{(i)}\left(\left(\frac{1}{r}\right)^-\right)$.
Then by the hypothesis of the theorem, we have $\be_i = 0$ and
$f^{(i)}\left(\left( \frac{1}{r} \right)^-\right)=0$ for all $i=0,1,2,\ldots,(\ell-1)$.
So the Taylor series expansions of
$f(1-z_n)$ around $z_n=0^+$ up to $\ell$ 
and $G(z_n)$ around $0^+$ up to order $(\ell+1)$
so that $z_n \in (0,\ve)$, are as follows:
$$f(1-z_n)=\frac{(-1)^{\ell}}{\ell!}f^{(\ell)}(1^-)\,z_n^{\ell}+O\left(z_n^{\ell+1}\right)$$
$$G(z_n)=G(0^+)+ \frac{1}{(\ell+1)!}\left(\frac{d^{\ell+1}G(0^+)}{d z_n^{\ell+1}}\right)\,z_n^{\ell+1}+
+O\left( z_n^{\ell+2} \right) \\
=1+\frac{(-1)^{\ell+1}\be_{\ell}}{(\ell+1)!}\,z_n^{\ell+1}+O\left(z_n^{\ell+2}\right).$$
Then substituting these expansions in Equation \eqref{eqn:asy-g=2-zn},
we get
$$
p_{{}_n}(F,r,1/r) \approx \int_0^{\ve}
n\Biggl[\frac{(-1)^{\ell}}{\ell!}f^{(\ell)}(1^-)\,z_n^{\ell}+O\left(z_n^{\ell+1}\right)\Biggr]
\Biggl[1-\frac{(-1)^{\ell}\be_{\ell}}{(\ell+1)!}\,z_n^{\ell+1}+O\left(z_n^{\ell+2}\right)\Biggr]^{n-1}\,d z_n.
$$

Now we let 
$z_n=v\,n^{-1/(\ell+1)}$,
to obtain
\begin{multline}
\label{eqn:Pg2-in-(0,1)-1/r}
p_{{}_n}(F,r,1/r) \approx
\int_0^{\ve\,n^{1/(\ell+1)}}n\,
\Biggl[\frac{(-1)^{\ell}}{n^{\ell/(\ell+1)}\,\ell!}f^{(\ell)}(1^-)v^{\ell}+O\left(n^{-1}\right)\Biggr]\\
\Biggl[1-\frac{1}{n}\left(\frac{(-1)^{\ell}\be_{\ell}}{(\ell+1)!}\,v^{\ell+1}+
O\left(n^{-(\ell+2)/(\ell+1)}\right)\right)\Biggr]^{n-1}\,\left(\frac{1}{n^{1/(\ell+1)}}\right)\,dv\\
\text{letting $n \rightarrow
\infty,$~~~~~~~~~~~~~~~~~~~~~~~~~~~~~~~~~~~~~~~~~~~~~~~~~~~~~~~~
~~~~~~~~~~~~~~~~~~~~~~~~~~~~~~~~~~~~~~~~~~~~~~~~~~~~}\\
\approx
\int_0^{\infty} \frac{(-1)^{\ell}}{\ell!} f^{(\ell)}(1^-)v^{\ell} \,
\exp\left[-\frac{(-1)^{\ell}\be_{\ell}}{(\ell+1)!}\,v^{\ell+1}\right]\, dv
= \frac{f^{(\ell)}(1^-)}{\be_{\ell}}
= \frac{f^{(\ell)}(1^-)}
{f^{(\ell)}(1^-)+r^{-(\ell+1)}\,f^{(\ell)} \left( \left( \frac{1}{r} \right)^- \right)},
\end{multline}
as $n \rightarrow \infty$ at rate $O(\kappa_2(f)\cdot n^{-(\ell+2)/(\ell+1)})$.

For the general case of $\Y=\{\y_1,\y_2\}$,
as in the proof of Theorem \ref{thm:kth-order-gen-(r-1)/r},
using the transformation
$\phi(x)=(x-\y_1)/(\y_2-\y_1)$ and replacing $f(x)$ by $g(x)$ in Equation \eqref{eqn:Pg2-in-(0,1)-1/r},
the desired result follows.
$\blacksquare$

\end{document}